\documentclass[11pt,english]{article}
\usepackage[T1]{fontenc}
\usepackage[latin9]{inputenc}
\usepackage{amsmath}
\usepackage{setspace}
\usepackage{amssymb}
\usepackage{esint}

\makeatletter
\date{}

\makeatother

\usepackage{babel}

\begin{document}

\title{Division Theorems for Exact Sequences}

\author{Qingchun Ji}
\maketitle
\begin{description}
\item [{Abstract.}] \noindent Under certain integrability and geometric
conditions, we prove division theorems for the exact sequences of
holomorphic vector bundles and improve the results in the case of
Koszul complex. By introducing a singular Hermitian structure on the
trivial bundle, our results recover Skoda's division theorem for holomorphic
functions on pseudoconvex domains in complex Euclidean spaces.
\item [{Mathematics$\;$Subject$\;$Classification(2010):}] \noindent 32A38,32W05,53C55.
\end{description}
\begin{center}
$\mathbf{Introduction}$
\par\end{center}

$\,$

A classical problem in complex and algebraic geometry is to characterize
the image of global holomorphic sections under a sheaf-homomorphism.
Let $E,E^{'}$ be holomorphic vector bundles over a complex manifold,
and $\Phi:E\rightarrow E^{'}$ be a holomorphic homomorphism, it is
interesting to characterize the image of the induced homomorphism
on cohomology groups

\begin{center}
$H^{q}(M,\Omega^{p}(E))\rightarrow H^{q}(M,\Omega^{p}(E^{'}))$
\par\end{center}

\noindent under effective integrability and differential-geometric
conditions. Skoda's Division theorem ({[}S72,S78,D82{]}) is one of
the fundamental result on this kind of questions. Skoda's theorem
has played important roles in many important work in algebraic geometry.
Siu used Skoda's theorem to establish the deformation invariance of
plurigenera and prove the finite generation of canonical ring of compact
complex algebraic manifolds of general type({[}Siu 98,00,04,05,07{]}).
Skoda's theorem is an analogue of Hilbert's Nullstellensatz, but the
remarkable feature of effectiveness makes it very powerful. In {[}B87,
Laz99{]}, the authors used Skoda's theorem to prove effective versions
of the Nullstellensatz.

The statement of Skoda's theorem is the following: Let $\Omega$ be a
pseudoconvex domain in\emph{ $\mathbb{C}^{n},\psi\in
\rm{PSH}(\Omega)$ }(the set of plurisubharmonic functions on
$\Omega$), $g_{1}\cdots,g_{r}\in\mathcal{O}(\Omega)$ (the set of
holomorphic functions on $\Omega$), then for every
$f\in\mathcal{O}(\Omega)$ with
$\int_{\Omega}|f|^{2}|g|^{-2(q+q\varepsilon+1)}e^{-\psi}dV<+\infty,$
there exist holomorphic functions
$h_{1},\cdots,h_{r}\in\mathcal{O}(\Omega)$ such that $f=\sum
g_{i}h_{i}$ holds on $\Omega$
and$\int_{\Omega}|h|^{2}|g|^{-2q(1+\varepsilon)}e^{-\psi}dV\leq\frac{1+\varepsilon}{\varepsilon}\int_{\Omega}|f|^{2}|g|^{-2(q+q\varepsilon+1)}e^{-\psi}dV$
where$\left|g\right|^{2}=\sum_{i}\left|g_{i}\right|^{2},\left|h\right|^{2}=\sum_{i}\left|h_{i}\right|^{2},q=\textrm{min}\{n,r-1\}$
and $\varepsilon>0$ is a constant. In {[}S78,D82{]}, this theorem
was generalized to (generic) surjective homomorphisms between
holomorphic vector bundles by solving
$\overline{\partial}$-equations. The surjectivity was used in two
places in {[}S78{]} and {[}D82{]}. One is to construct a smooth
splitting map which reduces the question under consideration to
$\overline{\partial}$- equations, another is to relate an arbitrary
Hermitian structure on the quotient bundle to the quotient Hermitian
structure induced by the surjective homomorphism(actually this was
also done by constructing a smooth splitting in a more general
case). A recent generalization of Skoda's theorem was given in
{[}V08{]} for line bundles where the notion of Skoda triple was
introduced, by choosing Skoda triples the author established various
division theorems.

In this paper, we consider the question of establishing division theorems
in more general settings. In section 4, we will prove division theorems
for a generically exact sequence of holomorphic vector bundles(theorem
4.2, corollaries 4.3, 4.4) and division theorems for a single homomorphism(theorem
4.2$^{'}$). Moreover, in section 5, for the case of the Koszul complex
which involves nonsurjective homomorphisms, we improve the result
for general exact sequences(theorem 5.4 and corollary 5.5).

As a corollary, we get a division theorem for Koszul complex over
a pseudoconvex domain in $\mathbb{C}^{n}.$ More precisely, let $\Omega$
be a domain in\emph{ $\mathbb{C}^{n},g_{1}\cdots,g_{r}\in\mathcal{O}(\Omega),$
}we define a sheaf-homomorphism for each $1\leq\ell\leq r$ by setting

\noindent \begin{center}
$\begin{array}{ccccc}
\wedge^{\ell}\mathcal{O}_{\Omega}^{\oplus r} & \rightarrow & \wedge^{\ell-1}\mathcal{O}_{\Omega}^{\oplus r}\\
(h_{i_{1}\cdots i_{\ell}})_{i_{1}\cdots i_{\ell}=1}^{r} & \mapsto & (f_{i_{1}\cdots i_{\ell-1}})_{i_{1}\cdots i_{\ell-1}=1}^{r} & \textrm{with}\ f_{i_{1}\cdots i_{\ell-1}}=\underset{1\leq\nu\leq r}{\sum}g_{\nu}h_{\nu i_{1}\cdots i_{\ell-1}}.\end{array}$\emph{ }
\par\end{center}

\begin{singlespace}
\noindent By choosing appropriate Hermitian structure and the standard
argument of weak compactness, we show that Skoda's theorem is exactly
the special case where $\ell=1$ of the following result(see corollary
5.7) about the homomorphisms defined in this way.
\end{singlespace}

\begin{singlespace}
Let $\Omega$ be a pseudoconvex domain in\emph{
$\mathbb{C}^{n},g_{1}\cdots,g_{r}\in\mathcal{O}(\Omega),\psi\in\rm{PSH}(\Omega)$
}and $\varepsilon>0$ a constant, then for every global section
$(f_{i_{1}\cdots i_{\ell-1}})_{i_{1}\cdots
i_{\ell-1}=1}^{r}\in\Gamma(\Omega,\wedge^{\ell-1}\mathcal{O}_{\Omega}^{\oplus
r})$ ($1\leq\ell\leq r$ ) satisfying
\end{singlespace}

\noindent \begin{center}
$\underset{1\leq\nu\leq r}{\sum}g_{\nu}f_{\nu i_{1}\cdots i_{\ell-2}}=0$
and $\int_{\Omega}|f|^{2}|g|^{-2(q+q\varepsilon+1)}e^{-\psi}dV<+\infty,$
\par\end{center}

\noindent there exists at least one $(h_{i_{1}\cdots i_{\ell}})_{i_{1}\cdots i_{\ell}=1}^{r}\in\Gamma(\Omega,\wedge^{\ell}\mathcal{O}_{\Omega}^{\oplus r})$
such that

\noindent \begin{center}
$f_{i_{1}\cdots i_{\ell-1}}=\underset{1\leq\nu\leq r}{\sum}g_{\nu}h_{\nu i_{1}\cdots i_{\ell-1}},$
\par\end{center}

\noindent and

\noindent \begin{center}
$\int_{\Omega}|h|^{2}|g|^{-2q(1+\varepsilon)}e^{-\psi}dV\leq\frac{1+\varepsilon}{\varepsilon}\int_{\Omega}|f|^{2}|g|^{-2(q+q\varepsilon+1)}e^{-\psi}dV,$
\par\end{center}

\noindent where $\left|g\right|^{2}=\sum_{i}\left|g_{i}\right|^{2},$
$\left|h\right|^{2}=\underset{i_{1}<\cdots<i_{\ell}}{\sum}\left|h_{i_{1}\cdots i_{\ell}}\right|^{2},$
$\left|f\right|^{2}=\underset{i_{1}<\cdots<i_{\ell-1}}{\sum}\left|f_{i_{1}\cdots i_{\ell-1}}\right|^{2},$
$q=\textrm{min}\{n,r-1\}.$

In the surjective case, we also give sufficient conditions for the
solvability of division problem in matrix form. This is achieved by
applying Skoda's original division theorem and a purely algebraic
argument(see corollary 5.9). Our integrability conditions are
comparable with those conditions appeared in {[}KT71{]}.

As an example of the application of Spencer's trick, we could construct
smooth lifting for an exact sequence and reduce the question to solving
$\overline{\partial}$- equations with additional restrictions, the
main difficulty would be that it requires to look for solutions of
$\overline{\partial}$-equations which take values in a subsheaf.
Although we don't use this fact essentially in our proof of division
theorems, we will sketch the rough ideal in section 4.

Instead of solving $\overline{\partial}$-equations, we hope to implement
Skoda's original estimate in general case. To this end, we first formulate
an algebraic inequality(lemma 2.2) which helps to complete square
in the proof of our main estimate(lemma 3.2). To apply this algebraic
lemma, we need to consider certain bundle homomorphisms, and the second
fundamental form is quite useful which is different from that used
in {[}S78{]} and {[}D82{]}. Since the bundle homomorphism is surjective
in {[}S78{]} and {[}D82{]}, the second fundamental form of the kernel
is well defined, but in more general case the kernel is only a subsheaf
which is no longer a subbundle. For this reason, we consider the second
fundamental form of the line bundle spanned by the homomorphism itself
inside the holomorphic bundle $\textrm{Hom}(E,E^{'})$ which is always
well defined outside the subset of zeros of the bundle-homomorphism
under consideration. To obtain our main estimate, we also use the
twisted Bochner-Kodaira-Nakano formula which is quite natural and
useful in geometry and analysis. Applications of this kind of technique
could be found in many references, e.g. {[}OT{]},{[}Siu82{]},{[}Siu00{]},{[}Siu04{]}
and {[}Dr08{]}. The discussion about the division problem of Koszul
complex via the residue theory and its interesting applications could
be found in {[}A04{]}, {[}A06{]} and {[}AG10{]}.

\section*{{\normalsize 1. The Second Fundamental Form}}

$\ \ \ $Let $(M,\omega)$ be a K\"{a}hler manifold,
$(E,h),(E^{'},h^{'})$ Hermitian holomorphic vector bundles over $M$,
$\textrm{dim}_{\mathbb{C}}M=n,\textrm{rank}_{\mathbb{C}}E=r,\textrm{rank}_{\mathbb{C}}E^{'}=r^{'}.$
$\Phi:\mathcal{O}_{M}(E)\rightarrow\mathcal{O}_{M}(E^{'})$ is a
sheaf-homomorphism, given by a holomorphic section
$\Phi\in\Gamma(M,\textrm{Hom}(E,E^{'})).$

The homomorphism $\Phi$ is assumed, for simplicity, to be nowhere
zero on $M$ in sections 1 and 2.

We shall adopt the following convention on the range of indices.

\begin{center}
$1\leq\alpha,\beta$$\leq n,1\leq i,j\leq r,1\leq a,b\leq r^{'}.$
\par\end{center}

\noindent We also adhere to the summation convention that sum is performed
over strictly increasing multi-indices.

By using local coordinates $z_{1},\cdots,z_{n}$ of \emph{$M$,} holomorphic
farmes\emph{ $\{e_{1},\cdots,e_{r}\},$ }$\{e_{1}^{'},\cdots,e_{r^{'}}^{'}\}$
of $E$ and $E^{'}$ respectively we write

\begin{center}
$\omega=\sqrt{-1}g_{\alpha\overline{\beta}}dz_{\alpha}\wedge d\overline{z}_{\beta},h=h_{i\overline{j}}e_{i}^{*}\otimes\overline{e_{j}^{*}},h^{'}=h_{a\overline{b}}^{'}e_{a}^{'*}\otimes\overline{e_{b}^{'*}},$
\par\end{center}

\noindent and

\begin{equation}
\Phi=\Phi_{ia}e_{a}^{'}\otimes e_{i}^{*}.\end{equation}

\noindent where \emph{$\{e_{1}^{*},\cdots,e_{r}^{*}\}$ }and\emph{
$\{e_{1}^{'*},\cdots,e_{r^{'}}^{'*}\}$} are dual frames of $E^{*}$
and $E^{'*}$ respectively. We also know by definition that the adjoint
homomorphism of $\Phi$ w.r.t. the Hermitian structures $h$ and $h^{'}$
is given by \begin{equation}
\Phi^{*}=\overline{\Phi_{ia}}h^{\overline{i}j}h_{a\overline{b}}^{'}e_{j}\otimes e_{b}^{'*}\end{equation}

\noindent where $(h^{\overline{i}j})=(h_{i\overline{j}})^{-1}.$

It is convenient to make it a convention that we always work with
normal coordinates and normal frames for a pointwise computation,
i.e. we have $g_{\alpha\overline{\beta}}=\delta_{\alpha\beta},h_{i\overline{j}}=\delta_{ij},h_{a\overline{b}}^{'}=\delta_{ab}$
and $dg_{\alpha\overline{\beta}}=dh_{i\overline{j}}=dh_{a\overline{b}}^{'}=0$
at the point under consideration.

We define a smooth section $B$ of \emph{$A^{1,0}({\rm
Hom}(E,E^{'}))$ }as the differentiation\emph{ }of \emph{$\Phi$,
}i.e. for every $X\in T^{1,0}M$\emph{ }\begin{equation}
BX=\nabla_{X}^{\textrm{Hom}(E,E^{'})}\Phi\end{equation}
 where $\nabla^{\textrm{Hom}(E,E^{'})}$ is the induced Chern connection
on $\textrm{Hom}(E,E^{'}).$

Let \emph{$B^{*}$ }be the adjoint of $B$ where $B$ is viewed as
a homomorphism from $T^{1,0}M$ to $\textrm{Hom}(E,E^{'}).$ By definition,
we know \emph{$B^{*}$ }is a smooth section of $\textrm{Hom}(\textrm{Hom}(E,E^{'}),T^{1,0}M).$
A pointwise computation shows that the exterior differentiation of
$\Phi^{*}$ and $B^{*}$are related as follows

\begin{equation}
\overline{\partial}\Phi^{*}=(B^{*}e_{i}^{*}\otimes e_{a}^{'})^{\flat}\otimes e_{i}\otimes e_{a}^{'*}\end{equation}
where $\flat$ is the \emph{$\mathbb{C}$-}linear map of lowering
indices by using the metric $g_{\alpha\overline{\beta}}$.

Now we compute the derivatives up to second order of the function
$\varphi=\log||\Phi||^{2}$ in terms of the homomorphism \emph{$B$.}
Since $\Phi$ is holomorphic, we have

\begin{spacing}{0.7}
$\partial_{\alpha}\varphi=\left\Vert \Phi\right\Vert ^{-2}(\nabla_{\alpha}^{\textrm{Hom}(E,E^{'})}\Phi,\Phi)$\begin{equation}
=e^{-\varphi}(B\frac{\partial}{\partial z_{\alpha}},\Phi)\ \ \ \ \ \ \ \ \ \ \ \ \ \ \ \ \ \ \ \ \ \ \ \ \ \ \ \ \ \ \ \ \ \ \ \ \ \ \ \ \ \ \ \ \ \ \ \ \ \ \ \ \ \ \ \,\end{equation}

\begin{equation}
\textrm{grad}^{0,1}\varphi=\left\Vert \Phi\right\Vert ^{-2}g^{\overline{\alpha}\beta}(B\frac{\partial}{\partial z_{\beta}},\Phi)\frac{\partial}{\partial\overline{z}_{\alpha}}\ \ \ \ \ \ \ \ \ \ \ \ \ \ \ \ \ \ \ \ \ \ \ \ \ \ \ \ \ \ \ \ \ \ \ \ \ \ \ \end{equation}

\end{spacing}

$\partial_{\alpha}\partial_{\overline{\beta}}\varphi=e^{-\varphi}[-\partial_{\overline{\beta}}\varphi(B\frac{\partial}{\partial z_{\alpha}},\Phi)+(\nabla_{\overline{\beta}}^{\textrm{Hom}(E,E^{'})}\nabla_{\alpha}^{\textrm{Hom}(E,E^{'})}\Phi,\Phi)$

$\ \ \ \ \ \ \ \ \ \ \ \ +(B\frac{\partial}{\partial z_{\alpha}},\nabla_{\beta}^{\textrm{Hom}(E,E^{'})}\Phi)]$

$\ \ \ \ \ \ \ \ \ =e^{-\varphi}[-\left\Vert \Phi\right\Vert ^{-2}\overline{(B\frac{\partial}{\partial z_{\beta}},\Phi)}(B\frac{\partial}{\partial z_{\alpha}},\Phi)+(B\frac{\partial}{\partial z_{\alpha}},B\frac{\partial}{\partial z_{\beta}})$

\begin{equation}
+(F_{\overline{\beta}\alpha}^{\textrm{Hom}(E,E^{'})}\Phi,\Phi)]\ \ \ \ \ \ \ \ \ \ \ \ \ \ \ \ \ \ \ \ \ \ \ \ \ \ \ \ \ \ \ \ \ \ \ \ \ \ \ \ \ \end{equation}

\noindent where $F_{\overline{\beta}\alpha}^{\textrm{Hom}(E,E^{'})}=\left[\nabla_{\overline{\beta}}^{\textrm{Hom}(E,E^{'})},\nabla_{\alpha}^{\textrm{Hom}(E,E^{'})}\right]$
is the curvature of the induced Chern connection on $\textrm{Hom}(E,E^{'}).$

Let \emph{P }be the orthonormal projection from $\textrm{Hom}(E,E^{'})$
onto the subbundle $\textrm{Span}_{\mathbb{C}}\{\Phi\}^{\bot}\subseteq\textrm{Hom}(E,E^{'}),$
we define

\begin{equation}
B_{\Phi}=P\circ B:T^{1,0}M\rightarrow\textrm{Span}_{\mathbb{C}}\{\Phi\}^{\perp}.\end{equation}

By definition $B_{\Phi}$ is the second fundamental form of the holomorphic
line bundle $\textrm{Span}_{\mathbb{C}}\{\Phi\}$ in $\textrm{Hom}(E,E^{'}).$

Let $L$ be a Hermitian holomorphic line bundle over $M,$
\{$\sigma\}$ be a holomorphic local frame of $L,$ we denote its dual
frame by $\{\sigma^{*}\}.$ For any\emph{
$v=v_{\overline{K}i}dz\wedge d\overline{z}_{K}\otimes\sigma\otimes
e_{i}\in A^{n,k}(L\otimes E)$, }we define the associated smooth
section $A_{v}$ of \emph{${\rm Hom}(\wedge^{n,k-1}TM\otimes
L^{*}\otimes E^{*},T^{1,0}M)$ }by setting

\begin{equation}
A_{v}(\frac{\partial}{\partial z}\otimes\frac{\partial}{\partial\overline{z}_{K}}\otimes\sigma^{*}\otimes e_{i}^{*})=(-1)^{n}g^{\overline{\alpha}\beta}v_{\overline{\alpha K}i}\frac{\partial}{\partial z_{\beta}},\end{equation}

\noindent for $1\leq k\leq n$, $|K|=k-1,$ here we denote

\begin{center}
$dz=dz_{1}\wedge\cdots\wedge dz_{n}\ \textrm{and}\ \frac{\partial}{\partial z}=\frac{\partial}{\partial z_{1}}\wedge\cdots\wedge\frac{\partial}{\partial z_{n}}.$
\par\end{center}

It is easy to see $A_{v}$ is well defined. Moreover, if $\left\{ z_{1},\cdots,z_{n}\right\} $,
\emph{$\left\{ e_{1},\cdots,e_{r}\right\} ,$} $\left\{ e_{1}^{'},\cdots,e_{r^{'}}^{'}\right\} $
and \{$\sigma\}$ are normal at the point $x\in M,$ then we have

\begin{spacing}{0.8}
\begin{equation}
(A_{v}(\frac{\partial}{\partial z}\otimes\frac{\partial}{\partial\overline{z}_{K}}\otimes\sigma^{*}\otimes e_{i}^{*}),X)=(v,X^{\flat}\wedge dz\wedge d\overline{z}_{K}\otimes\sigma\otimes e_{i})\end{equation}

\end{spacing}

\noindent for every $X\in T_{x}^{1,0}M$ and multi-index $K$ with
$|K|=k-1$.

Multiplying both sides of the equality (7) by $v_{\overline{\alpha K}i}\overline{v_{\overline{\beta K}i}}$
and summing over $\alpha,\beta,i$ and increasing multi-indices $K$
with $|K|=k-1$ give the following expression of $\partial_{\alpha}\partial_{\overline{\beta}}\varphi v_{\overline{\alpha K}i}\overline{v_{\overline{\beta K}i}}$
which will be used to handle the curvature term in the Bochner-Kodaira-Nakano
formula.

$\partial_{\alpha}\partial_{\overline{\beta}}\varphi v_{\overline{\alpha K}i}\overline{v_{\overline{\beta K}i}}=e^{-\varphi}[-|(BA_{v}(\frac{\partial}{\partial z}\otimes\frac{\partial}{\partial\overline{z}_{K}}\otimes\sigma^{*}\otimes e_{i}^{*}),\frac{\Phi}{\left\Vert \Phi\right\Vert })|^{2}$

$\ \ \ \ \ \ \ \ \ \ \ \ \ \ \ \ \ \ \ \ \ \ \ \ +\left\Vert BA_{v}(\frac{\partial}{\partial z}\otimes\frac{\partial}{\partial\overline{z}_{K}}\otimes\sigma^{*}\otimes e_{i}^{*})\right\Vert ^{2}$

$\ \ \ \ \ \ \ \ \ \ \ \ \ \ \ \ \ \ \ \ \ \ \ \ +(F_{\overline{A_{v}(\frac{\partial}{\partial z}\otimes\frac{\partial}{\partial\overline{z}_{K}}\otimes\sigma^{*}\otimes e_{i}^{*})}A_{v}(\frac{\partial}{\partial z}\otimes\frac{\partial}{\partial\overline{z}_{K}}\otimes\sigma^{*}\otimes e_{i}^{*})}^{\textrm{Hom}(E,E^{'})}\Phi,\Phi)]$

$\ \ \ \ \ \ \ \ \ \ \ \ \ \ \ \ \ \ \ \ \ \ =e^{-\varphi}\left\Vert B_{\Phi}A_{v}\right\Vert ^{2}$

\begin{equation}
\ \ \ \ \ \ \ \ \ \ \ \ \ \ \ \ \ \ \ \ \ \ \ \ \ \ \ \ \ \ \ \ \ \ -e^{-\varphi}(F_{A_{v}(\frac{\partial}{\partial z}\otimes\frac{\partial}{\partial\overline{z}_{K}}\otimes\sigma^{*}\otimes e_{i}^{*})\overline{A_{v}(\frac{\partial}{\partial z}\otimes\frac{\partial}{\partial\overline{z}_{K}}\otimes\sigma^{*}\otimes e_{i}^{*})}}^{\textrm{Hom}(E,E^{'})}\Phi,\Phi).\end{equation}

\section*{{\normalsize 2. Algebraic Preliminaries}}

$\ \ \ $In this section we introduce the notion of trace which
generalizes the usual conception of trace for linear
transformations. By using the Cauchy-Schwarz inequality, we will
establish a fundamental estimate for the generalized trace which
plays an important role in our main estimate. In order to apply this
inequality, we also collect in this section the pointwise
calculations concerning the smooth section
$\textrm{Tr}B_{\Phi}A\in\Gamma(\wedge_{M}^{n,k-1}\otimes L\otimes
E^{'}).$ Lemma 2.2 and lemma 2.3 make up the major algebraic part of
the proof of our main estimate. Its geometric ingredient is the
twisted Bochner-Kodaira-Nakano formula(combined with the Morrey's
trick) on a bounded domain which will be discussed in section 3.\\

\noindent $\mathbf{Definition2.1.}$ Let $U,V,W$ be linear space\emph{s},
$D\in\textrm{Hom}(U,V)$, $\rho:V\times U^{*}\rightarrow W$ a bilinear
map. If $U$ is finite-dimensional we define the trace $\textrm{Tr}_{\rho}D\in W$
of the linear map $D$ w.r.t. $\rho$ to be

\begin{equation}
\textrm{Tr}_{\rho}D=\underset{i}{\sum}\rho(Du_{i},u^{i})\end{equation}

\noindent where \emph{$\{u_{i}\}$ }is a basis of \emph{$U$,} \emph{$\{u^{i}\}\subseteq U^{*}$
}is its dual basis.

$\,$

The definition of $Tr_{\rho}D$ is obviously independent of the
choice of the basis \emph{$\{u_{i}\}.$ }Now we establish the basic
estimate for $Tr_{\rho}D$.\\

\noindent $\mathbf{Lemma2.2.}$ If $U,V,W$ are Hermitian space\emph{s},
$D\in\textrm{Hom}(U,V)$ and $\rho:V\times U^{*}\rightarrow W$ is
a bilinear map, then we have

\begin{equation}
\left\Vert \textrm{Tr}_{\rho}D\right\Vert _{W}\leq\sqrt{\textrm{rank}(D)}\left\Vert \rho\right\Vert \left\Vert D\right\Vert .\end{equation}

$\,$

\noindent Proof. By using the singular value decomposition theorem
for a linear map between Hermitian spaces, one can always find an
orthonormal basis \emph{$\{u_{i}\}$} of\emph{ $U$} such that

\noindent \begin{center}
$(Du_{i},Du_{j})_{V}=0$ for $i\neq j$.
\par\end{center}

\noindent Since $\textrm{Im}(D)=\textrm{Span}_{\mathbb{C}}\{Du_{i}\}$
and $Du_{1},Du_{2},\cdots$ are mutually perpendicular, we have

\noindent \begin{center}
$\sharp\{i\mid Du_{i}\neq0\}=\textrm{dim(Im}(D))=\textrm{rank}(D).$
\par\end{center}

\noindent Let \emph{$\{u^{i}\}\subseteq U^{*}$ }be the dual basis
of \emph{$\{u_{i}\}$}, then $\{u^{i}\}$ forms a orthonormal basis
of \emph{$U^{*},\textrm{i.e.}(u^{i},u^{j})=\delta_{ij}.$} Consequently,
by using the Cauchy-Schwarz inequality and the definition of $\left\Vert \rho\right\Vert $:

\noindent \begin{center}
$\left\Vert \rho\right\Vert :=\underset{\left\Vert \alpha\right\Vert _{U^{*}}=1}{\underset{\left\Vert v\right\Vert _{V}=1}{\textrm{max}}}\left\Vert \rho(v,\alpha)\right\Vert _{W}\:\textrm{and}\:\left\Vert D\right\Vert =\sqrt{\underset{i}{\sum}\left\Vert Du_{i}\right\Vert _{V}^{2}},$
\par\end{center}

\noindent we get the desired estimate as follows

\begin{center}
$\left\Vert \textrm{Tr}_{\rho}D\right\Vert _{W}^{2}=\left\Vert \underset{i}{\sum}\rho(Du_{i},u^{i})\right\Vert _{W}^{2}$
\par\end{center}

\begin{center}
$\ \ \ \ \ \ \ \ \ \ \ \ \ \ \ \ \ \ \leq\underset{i}{(\sum}\left\Vert \rho(Du_{i},u^{i})\right\Vert _{W})^{2}$
\par\end{center}

\begin{center}
$\ \ \ \ \ \ \ \ \ \ \ \ \ \ \ \ \leq\left\Vert \rho\right\Vert ^{2}\underset{i}{(\sum}\left\Vert Du_{i}\right\Vert _{V})^{2}$
\par\end{center}

\begin{center}
$\ \ \ \ \ \ \ \ \ \ \ \ \ \ \ \ \ \ \ \ =\left\Vert \rho\right\Vert ^{2}(\underset{Du_{i}\neq0}{\sum}\left\Vert Du_{i}\right\Vert _{V})^{2}$
\par\end{center}

\begin{center}
$\ \ \ \ \ \ \ \ \ \ \ \ \ \ \:\ \ \ \ \ \ \ \ \ \ \ \ \ \leq\textrm{rank}(D)\left\Vert \rho\right\Vert ^{2}\underset{Du_{i}\neq0}{\sum}\left\Vert Du_{i}\right\Vert _{V}^{2}$
\par\end{center}

\begin{center}
$\ \ \ \ \ \ \ \:\ \ \ \ \ \ \ \ \ \ =\textrm{rank}(D)\left\Vert \rho\right\Vert ^{2}\left\Vert D\right\Vert ^{2},$
\par\end{center}

\begin{spacing}{0.7}
\noindent which completes the proof.
\end{spacing}

\begin{spacing}{0.7}
\noindent \begin{flushright}
$\square$
\par\end{flushright}
\end{spacing}

We will apply lemma 2.2 in two specific circumstances. In the following
lemma 2.3, we will choose $V=W\otimes U$ and $\rho:W\otimes U\times U^{*}\rightarrow W$
to be the natural contraction between $U$ and its dual space $U{}^{*}.$
If we identify $U$ and $U^{*}$ via the $\mathbb{C}$-antilinear
isomorphism defined by the Hermitian structure on $U,$ then we have
explicitly

\begin{center}
$\textrm{Tr}_{\rho}D=(Du_{i},w_{a}\otimes u_{i})_{W\otimes U}w_{a}$
\par\end{center}

\noindent for orthonormal bases $\{u_{i}\}\subseteq U,\{w_{a}\}\subseteq W.$

In section 5, we will consider $U=\wedge^{p}V,$ and $\rho:V\times\wedge^{p}V^{*}\rightarrow\wedge^{p-1}V^{*}$defined
by the interior product, i.e.

\begin{center}
$\rho(v,\xi):=v\lrcorner\xi.$
\par\end{center}

Obviously, the Cauchy-Schwarz inequality shows that $\left\Vert
\rho(v,u^{*})\right\Vert _{W}\leq\left\Vert v\right\Vert
_{V}\left\Vert u^{*}\right\Vert _{U^{*}},$ we get therefore
$\left\Vert \rho\right\Vert \leq1$ in both of the cases mentioned
above. Since we always work with specific bilinear map, the
subscript $\rho$ will be omitted without causing ambiguity. The
importance of the inequality (13) is that the coefficient of the
left hand side only depends on the rank of $D.$\\

Now we proceed to prove the key identity involving $\textrm{Tr}B_{\Phi}A.$
Since the computations in this section are pointwise, as mentioned
before, we will work with normal coordinates and normal frames at
a given point.

For a given $v\in A^{n,k}(L\otimes E),1\leq k\leq n,$ we define by
(9) the associated homomorphism \emph{$A_{v}\in{\rm
Hom}(\wedge^{n,k-1}TM\otimes L^{*}\otimes E^{*},T^{1,0}M)$. }Under
the standard bundle isomorphism

\begin{singlespace}
\begin{center}
$\textrm{Hom}(\wedge^{n,k-1}TM\otimes L^{*}\otimes E^{*},\textrm{Hom}(E,E^{'}))$\begin{equation}
\cong\textrm{Hom}(E^{*},\wedge^{n,k-1}M\otimes L\otimes E^{'}\otimes E^{*}),\ \ \ \ \ \ \ \ \ \ \end{equation}

\par\end{center}
\end{singlespace}

\noindent we could define $\textrm{Tr}B_{\Phi}A_{v}\in\Gamma(\wedge^{n,k-1}M\otimes L\otimes E^{'})$
by (12) where the bilinear map $\rho$ is given by the pairing between
$E$ and $E^{*}$ which is defined by the Hermitian structure on $E$.

$\,$

The main result of about $\textrm{Tr}B_{\Phi}A_{v}$ is recorded in
the following formula.\\

\noindent $\mathbf{Lemma2.3.}$ For any $u$ in $\wedge^{n,k-1}M\otimes L\otimes E^{'}$
and $v$ in $\wedge^{n,k}M\otimes L\otimes E,$ we have \begin{equation}
(\overline{\partial}\Phi^{*}\wedge u,v)-(\Phi^{*}u,\textrm{grad}^{0,1}\varphi\lrcorner v)=(u,\textrm{Tr}B_{\Phi}A_{v})\end{equation}

\noindent where $\varphi=\log\left\Vert \Phi\right\Vert ,\Phi\in\Gamma(M,\textrm{H\textrm{om}}(E,E^{'}))$
and $A_{v}$ is defined by $v$ as described in (9), $1\leq k\leq n$.

$\ $

\noindent Proof. Let $u=u_{\overline{K}a}dz\otimes d\overline{z}_{K}\otimes\sigma\otimes e_{a}^{'}\in\wedge^{n,k-1}M\otimes L\otimes E^{'},$
$v=v_{\overline{J}i}dz\otimes d\overline{z}_{J}\otimes\sigma\otimes e_{i}\in\wedge^{n,k}M\otimes L\otimes E.$
We know by the definition (12) and the identification (14) that

$\textrm{Tr}B_{\Phi}A_{v}=(e_{i}^{*}\lrcorner B_{\Phi}A_{v},dz\otimes d\overline{z}_{K}\otimes\sigma\otimes e_{a}^{'}\otimes e_{i}^{*})dz\otimes d\overline{z}_{K}\otimes\sigma\otimes e_{a}^{'}$

$\,\ \ \ \ \ \ \ \ \ \ \ \:=(B_{\Phi}A_{v}(\frac{\partial}{\partial z}\otimes\frac{\partial}{\partial\overline{z}_{K}}\otimes\sigma^{*}\otimes e_{i}^{*}),e_{a}^{'}\otimes e_{i}^{*})dz\otimes d\overline{z}_{K}\otimes\sigma\otimes e_{a}^{'},$

\noindent which implies that

\begin{singlespace}
$(u,\textrm{Tr}B_{\Phi}A_{v})=u_{\overline{K}a}(dz\otimes d\overline{z}_{K}\otimes\sigma\otimes e_{a}^{'},\textrm{Tr}B_{\Phi}A_{v})$

\begin{equation}
\,\ \ \ \ =u_{\overline{K}a}(e_{a}^{'}\otimes e_{i}^{*},B_{\Phi}A_{v}(\frac{\partial}{\partial z}\otimes\frac{\partial}{\partial\overline{z}_{K}}\otimes\sigma^{*}\otimes e_{i}^{*})).\end{equation}

\end{singlespace}

\noindent From the equalities (4) and (10), it follows that

$(\overline{\partial}\Phi^{*}\wedge u,v)=((B^{*}(e_{a}^{'}\otimes e_{i}^{*}))^{\flat}\otimes e_{i}\otimes e_{a}^{'*}(u),v)$

$\ \ \ \ \ \ \ \ \ \ \ \ \ \ \ \ \,=u_{\overline{K}b}((B^{*}(e_{a}^{'}\otimes e_{i}^{*}))^{\flat}\wedge dz\otimes d\overline{z}_{K}\otimes\sigma\otimes e_{a}^{'*}(e_{b}^{'})e_{i},v)$

$\ \ \ \ \ \ \ \ \ \ \ \ \ \ \ \ \,=u_{\overline{K}a}((B^{*}(e_{a}^{'}\otimes e_{i}^{*}))^{\flat}\wedge dz\otimes d\overline{z}_{K}\otimes\sigma\otimes e_{i},v)$

$\ \ \ \ \ \ \ \ \ \ \ \ \ \ \ \ \,=u_{\overline{K}a}(B^{*}(e_{a}^{'}\otimes e_{i}^{*}),A_{v}(\frac{\partial}{\partial z}\otimes\frac{\partial}{\partial\overline{z}_{K}}\otimes\sigma^{*}\otimes e_{i}^{*}))$

\begin{equation}
\,\ =u_{\overline{K}a}(e_{a}^{'}\otimes e_{i}^{*},BA_{v}(\frac{\partial}{\partial z}\otimes\frac{\partial}{\partial\overline{z}_{K}}\otimes\sigma^{*}\otimes e_{i}^{*})).\end{equation}

\noindent We also obtain form (6) that

$\ \ \ (\Phi^{*}u,\textrm{grad}^{0,1}\varphi\lrcorner v)$

$=(\Phi^{*}u,\left\Vert \Phi\right\Vert ^{-2}(B\frac{\partial}{\partial z_{\alpha}},\Phi)\frac{\partial}{\partial\overline{z}_{\alpha}}\lrcorner v)$

$=(-1)^{n}(\Phi^{*}u,\left\Vert \Phi\right\Vert ^{-2}(B\frac{\partial}{\partial z_{\alpha}},\Phi)v_{\overline{\alpha K}i}dz\otimes d\overline{z}_{K}\otimes\sigma\otimes e_{i})$

$=(\Phi^{*}u,\left\Vert \Phi\right\Vert ^{-2}(BA_{v}(\frac{\partial}{\partial z}\otimes\frac{\partial}{\partial\overline{z}_{K}}\otimes\sigma^{*}\otimes e_{i}^{*}),\Phi)dz\otimes d\overline{z}_{K}\otimes\sigma\otimes e_{i})$

$=(u,(BA_{v}(\frac{\partial}{\partial z}\otimes\frac{\partial}{\partial\overline{z}_{K}}\otimes\sigma^{*}\otimes e_{i}^{*}),\frac{\Phi}{\left|\Phi\right|})dz\otimes d\overline{z}_{K}\sigma\otimes\frac{\Phi}{\left|\Phi\right|}(e_{i}))$

$=(u,\frac{\Phi_{ib}}{\left|\Phi\right|}(BA_{v}(\frac{\partial}{\partial z}\otimes\frac{\partial}{\partial\overline{z}_{K}}\otimes\sigma^{*}\otimes e_{i}^{*}),\frac{\Phi}{\left|\Phi\right|})dz\otimes d\overline{z}_{K}\otimes\sigma\otimes e_{b}^{'})$

\begin{spacing}{0.7}
$=u_{\overline{K}a}\frac{\Phi_{ia}}{\left|\Phi\right|}(\frac{\Phi}{\left|\Phi\right|},BA_{v}(\frac{\partial}{\partial z}\otimes\frac{\partial}{\partial\overline{z}_{K}}\otimes\sigma^{*}\otimes e_{i}^{*}))$

\begin{equation}
=u_{\overline{K}a}(e_{a}^{'}\otimes e_{i}^{*},QBA_{v}(\frac{\partial}{\partial z}\otimes\frac{\partial}{\partial\overline{z}_{K}}\otimes\sigma^{*}\otimes e_{i}^{*}))\ \ \ \ \ \ \ \ \ \ \ \ \ \ \ \ \ \ \ \ \ \ \ \ \ \end{equation}

\end{spacing}

\begin{spacing}{0.7}
\noindent where \emph{Q }is the orthonormal projection from $E^{'}\otimes E^{*}$
onto the line bundle $\textrm{Span}_{\mathbb{C}}\{\Phi\}.$ Now (15)
follows from (16) (17) (18) and the definition (8) of $B_{\Phi}$,
the proof is complete.
\end{spacing}

\begin{spacing}{0.7}
\noindent \begin{flushright}
$\ $$\square$
\par\end{flushright}
\end{spacing}

\section*{{\normalsize 3. The Main Estimate}}

$\mathbf{Notations.}$$\ \ \ $ We introduce some notations which are
needed to simplify our statements. Given a measurable function $\mu$
on a K\"{a}hler manifold $(M,\omega),$ we define the associated
signed measure by setting

\begin{equation}
dV_{\mu}=\mu dV_{\omega}\end{equation}

\noindent where $dV_{\omega}=\frac{\omega^{n}}{n!}$ is the volume
form of $\omega.$ Let $\Omega$ be a domain in $M$. We denote by

\noindent \begin{center}
$(\cdot,\cdot)_{\Omega,\mu},\ ||\cdot||_{\Omega,\mu}$
\par\end{center}

\noindent the $L^{2}$-inner product and $L^{2}$-norm defined by
using the measure $dV_{\mu}$. When $\mu$ is nonnegative, the corresponding
Hilbert space of square integrable $(n,k)$-forms on $\Omega$ valued
in $L\otimes E$ and $L\otimes E^{'}$ will be denoted respectively
by $L_{n,k}^{2}(\Omega,L\otimes E,dV_{\mu})$ and $L_{n,k}^{2}(\Omega,L\otimes E^{'},dV_{\mu})$
respectively. The subscript $\mu$ in $(\cdot,\cdot){}_{\Omega,\mu}$
will be dropped for $\mu=1.$

$\:$

We recall the definition of $m$-tensor positivity.\\

\noindent $\mathbf{Definition3.1.}$ A Hermitian holomorphic vector
bundle $(E,h)$ is said to be $m$-tensor semi-positive(semi-negative)
if the curvature $F$ (of Chern connection ) satisfies $\sqrt{-1}F(\eta,\eta)\geq0(\leq0)$
for every $\eta=\eta_{\alpha i}\frac{\partial}{\partial z_{\alpha}}\otimes e_{i}\in T^{1,0}M\otimes E$
with $\textrm{rank}(\eta_{\alpha i})\leq m$ where $z_{1},\cdots,z_{n}$
are holomorphic coordinates of \emph{$M$, $\{e_{1},\cdots,e_{r}\}$}
is a holomorphic frame of $E$ and $m$ is a positive integer. In
this case, we write $E\geq_{m}0(E\leq_{m}0).$

$\:$

It is easy to see the above definition does not depend on the choice
of the holomorphic coordinates $z_{1},\cdots,z_{n}$ or the
holomorphic frame $\{e_{1},\cdots,e_{r}\}.$\\

\noindent $\mathbf{Definition3.2.}$ Let $E$ be a holomorphic vector
bundle over $M,$ $Z\subsetneqq M$ be a subvariety, and $h$ be a
Hermitian structure on $E|_{M\setminus Z}.$ If for each $z\in Z$,
there exist a neighborhood $U$ of $z$, a smooth frame
$\{e_{1},\cdots,e_{r}\}$ over $U$ and some constant $\kappa>0$ such
that the matrix
$\left[h_{i\overline{j}}(w)-\kappa\delta_{ij}\right]$ is
semi-positive for every $w\in U\setminus Z$ where
$h_{i\overline{j}}:=h(e_{i},e_{j})$ and $\delta_{ij}$ is the
Kronecker delta, then we call $h$ a singular Hermitian structure on
$E$ which has singularities in $Z.$

$\,$

Let $\Omega\Subset M$ be a domain with smooth boundary, $\rho\in C^{\infty}(\overline{\Omega})$
a defining function of $\Omega$, i.e. $\Omega$ is given by $\rho<0$
and $d\rho\neq0$ on $\partial\Omega$. $\Omega$ is said to be pseudoconvex
if the levi form $L_{\rho}$ is semi-positive on $T^{1,0}\partial\Omega$.
This condition is independent of the choice of the defining function.

$\:$

Now we are in the position to prove the main estimate:\\

\noindent $\mathbf{Lemma3.3.}$ Let $(M,\omega)$ be a K\"{a}hler
manifold, and let $E$ be a Hermitian holomorphic vector bundle over
$M$, $L$ a Hermitian holomorphic line bundle over $M$. The Hermitian
structures of these bundles may have singularity along
$\Phi^{-1}(0)$ and $\Omega\Subset M\setminus\Phi^{-1}(0)$ is a
pseudoconvex domain with smooth boundary. Assume that the following
conditions hold on $\Omega:$

1. $E\geq_{m}0,m\geq\textrm{min}\{n-k+1,r\},1\leq k\leq n$;

2. the curvature of $\textrm{Hom}(E,E^{'})$ satisfies

\begin{center}
$(F_{X\overline{X}}^{\textrm{Hom}(E,E^{'})}\Phi,\Phi)\leq0$ for every
$X\in T^{1,0}M$;
\par\end{center}

3. the curvature of $L$ satisfies

\begin{center}
$\sqrt{-1}(\varsigma c(L)-\partial\overline{\partial}\varsigma-\tau^{-1}\partial\varsigma\wedge\overline{\partial}\varsigma)\geq\sqrt{-1}q(\varsigma+\delta)\partial\overline{\partial}\varphi$.
\par\end{center}

\noindent Then the following estimate \begin{equation}
\left\Vert |\Phi|^{-2}\Phi^{*}u+\overline{\partial}^{*}v\right\Vert _{\Omega,\varsigma+\tau}^{2}+\left\Vert \overline{\partial}v\right\Vert _{\Omega,\varsigma}^{2}\geq\left\Vert u\right\Vert _{\Omega,\frac{\varsigma(\lambda\delta+\lambda\varsigma-\varsigma)}{(\varsigma+\delta)|\Phi|^{2}}}^{2}\end{equation}
holds for every $\overline{\partial}$-closed $u\in A^{n,k-1}(\overline{\Omega},L\otimes E)$
satisfying $|\Phi^{*}u|^{2}\geq\lambda|\Phi|^{2}|u|^{2}$ a.e.(w.r.t.$dV_{\omega}$)
on $\Omega$ and every $v\in A^{n,k}(\overline{\Omega},L\otimes E)\cap\textrm{Dom}(\overline{\partial}^{*}),$
where $c(L)$ denotes the Chern form, $q=\underset{\Omega}{\textrm{max}}\:\textrm{rank}B_{\Phi},\varphi=\log|\Phi|^{2},$
$0<\varsigma\in C^{\infty}(\overline{\Omega})$ and $\lambda,\delta,\tau$
are measurable functions on $\Omega$ satisfying $\lambda,\tau>0,\varsigma+\delta\geq0.$
All the weighted norms are described at the beginning of this section.

$\ $

\noindent Proof. Since we work on a fixed domain, the subscript $\Omega$
will be omitted in the following proof. We assume $\Omega$ is given
by $\rho<0$ and $\ d\rho\neq0$ on $\partial\Omega$ where $\rho\in C^{\infty}(\overline{\Omega})$.

Step 1. From the condition

\noindent \begin{center}
$|\Phi^{*}u|^{2}\geq\lambda|\Phi|^{2}|u|^{2},$
\par\end{center}

\noindent we get

$\,\,\ \ \ \ \ \ \ \ \ \ \ \ \ \ \left\Vert |\Phi|^{-2}\Phi^{*}u+\overline{\partial}^{*}v\right\Vert _{\Omega,\varsigma+\tau}^{2}+\left\Vert \overline{\partial}v\right\Vert _{\Omega,\varsigma}^{2}$

$\ \ \ \ \ \ \ \ \ \ \ \ =\left\Vert \sqrt{\varsigma+\tau}|\Phi|^{-2}\Phi^{*}u\right\Vert ^{2}+2\textrm{Re}((\varsigma+\tau)|\Phi|^{-2}\Phi^{*}u,\overline{\partial}^{*}v)$

$\ \ \ \ \ \ \ \ \ \ \ \ \ \ \ +\left\Vert \sqrt{\varsigma+\tau}\overline{\partial}^{*}v\right\Vert ^{2}+\left\Vert \sqrt{\varsigma}\overline{\partial}v\right\Vert ^{2}$

$\ \ \ \ \ \ \ \ \ \ \ \ =\left\Vert \sqrt{\varsigma}|\Phi|^{-2}\Phi^{*}u\right\Vert ^{2}+\left\Vert |\Phi|^{-2}\Phi^{*}u+\overline{\partial}^{*}v\right\Vert _{\tau}^{2}$

$\ \ \ \ \ \ \ \ \ \ \ \ \ \ \ +2\textrm{Re}(\varsigma e^{-\varphi}\Phi^{*}u,\overline{\partial}^{*}v)+\left\Vert \sqrt{\varsigma}\overline{\partial}^{*}v\right\Vert ^{2}+\left\Vert \sqrt{\varsigma}\overline{\partial}v\right\Vert ^{2}$

$\:$$\ \ \ \ \ \ \ \ \ \ \ \geq\left\Vert u\right\Vert _{|\Phi|^{-2}\lambda\varsigma}^{2}+\left\Vert |\Phi|^{-2}\Phi^{*}u+\overline{\partial}^{*}v\right\Vert _{\tau}^{2}$

\begin{equation}
+2\textrm{Re}(\varsigma e^{-\varphi}\Phi^{*}u,\overline{\partial}^{*}v)+\left\Vert \sqrt{\varsigma}\overline{\partial}^{*}v\right\Vert ^{2}+\left\Vert \sqrt{\varsigma}\overline{\partial}v\right\Vert ^{2}.\ \end{equation}

\noindent From the condition that

\noindent \begin{center}
$v\in A^{n,k}(\overline{\Omega},L\otimes E)\cap\textrm{Dom}(\overline{\partial}^{*}),$
\par\end{center}

\noindent we know

\noindent \begin{center}
$\textrm{grad}^{0,1}\rho\lrcorner v=0\ \textrm{on}\ \partial\Omega.$
\par\end{center}

\noindent By using the twisted Bochner-Kodaira-Nakano formula and
Morrey's trick, it follows from the pseudoconvexity of $\Omega$ that(see
{[}Siu82{]}, {[}Siu00{]})

$\ \ \ \left\Vert \sqrt{\varsigma}\overline{\partial}^{*}v\right\Vert ^{2}+\left\Vert \sqrt{\varsigma}\overline{\partial}v\right\Vert ^{2}$

$=\int_{\Omega}\varsigma|\overline{\nabla}v|^{2}+\varsigma(\left[\sqrt{-1}F^{L\otimes E},\Lambda_{\omega}\right]v,v)-\nabla^{\overline{\alpha}}\nabla^{\beta}\varsigma(\frac{\partial}{\partial\overline{z}_{\alpha}}\lrcorner v,\frac{\partial}{\partial\overline{z}_{\beta}}\lrcorner v)$

$\ \ \ +2\textrm{Re}(\textrm{grad}^{0,1}\varsigma\lrcorner v,\overline{\partial}^{*}v)dV_{\omega}+\int_{\partial\Omega}\varsigma\nabla^{\overline{\alpha}}\nabla^{\beta}\rho(\frac{\partial}{\partial\overline{z}_{\alpha}}\lrcorner v,\frac{\partial}{\partial\overline{z}_{\beta}}\lrcorner v)$

$\ \ \ +(\overline{\partial}\rho\wedge\textrm{grad}^{0,1}\varsigma\lrcorner v-\varsigma\overline{\partial}\rho\wedge\overline{\partial}^{*}v-\varsigma\overline{\partial}(\textrm{grad}^{0,1}\rho\lrcorner v),v)dA$

$=\int_{\Omega}\varsigma|\overline{\nabla}v|^{2}+\varsigma(\left[\sqrt{-1}F^{L\otimes E},\Lambda_{\omega}\right]v,v)-\nabla^{\overline{\alpha}}\nabla^{\beta}\varsigma(\frac{\partial}{\partial\overline{z}_{\alpha}}\lrcorner v,\frac{\partial}{\partial\overline{z}_{\beta}}\lrcorner v)$

$\ \ \ +2\textrm{Re}(\textrm{grad}^{0,1}\varsigma\lrcorner v,\overline{\partial}^{*}v)dV_{\omega}+\int_{\partial\Omega}\varsigma\nabla^{\overline{\alpha}}\nabla^{\beta}\rho(\frac{\partial}{\partial\overline{z}_{\alpha}}\lrcorner v,\frac{\partial}{\partial\overline{z}_{\beta}}\lrcorner v)dA$

$\geq\int_{\Omega}\varsigma(\left[\sqrt{-1}F^{L\otimes E},\Lambda_{\omega}\right]v,v)-\nabla^{\overline{\alpha}}\nabla^{\beta}\varsigma(\frac{\partial}{\partial\overline{z}_{\alpha}}\lrcorner v,\frac{\partial}{\partial\overline{z}_{\beta}}\lrcorner v)$

\begin{equation}
+2\textrm{Re}(\textrm{grad}^{0,1}\varsigma\lrcorner v,\overline{\partial}^{*}v)dV_{\omega},\ \ \ \ \ \ \ \ \ \ \ \ \ \ \ \ \ \ \ \ \ \ \ \ \ \ \ \ \ \ \ \ \ \ \ \ \ \ \ \ \ \ \ \end{equation}

\noindent $ $where $F^{L\otimes E}$ is the curvature of the Chern
connection on $L\otimes E,$ $\Lambda_{\omega}$ is the dual Lefschetz
operator of the K\"{a}hler form $\omega,$ and $dA$ is the induced
volume form on $\partial\Omega.$

$\,$Step 2. In order to obtain pointwise the lower bound of the integrand
of (22), we fix a point $x\in M$ and choose the local frames $\{e_{1},\cdots,e_{r}\},$
$\{\sigma\}$ of $E$ and $L$ respectively such that $(e_{i},e_{j})=\delta_{ij},|\sigma|^{2}=1$
at $x$. The following pointwise computations are carried out at this
fixed point $x.$

\noindent Set $F_{\alpha\overline{\beta}i\overline{j}}^{L\otimes E}=(F_{\alpha\overline{\beta}}^{L\otimes E}\sigma\otimes e_{i},\sigma\otimes e_{j})$,
then we have

$(\left[\sqrt{-1}F^{L\otimes E},\Lambda_{\omega}\right]v,v)=F_{\alpha\overline{\beta}i\overline{j}}^{L\otimes E}v_{\overline{\alpha K}i}\overline{v_{\overline{\beta K}j}}$

$\ \ \ \ \ \ \ \ \ \ \ \ \ \ \ \ \ \ \ \ \ \ \ \ \ \ \ \ \ \ =F_{\alpha\overline{\beta}i\overline{j}}^{E}v_{\overline{\alpha K}i}\overline{v_{\overline{\beta K}j}}+F_{\alpha\overline{\beta}}^{L}v_{\overline{\alpha K}i}\overline{v_{\overline{\beta K}i}}$

\noindent where $v=v_{\overline{J}i}dz\wedge d\overline{z}_{J}\otimes\sigma\otimes e_{i},$
$F_{\alpha\overline{\beta}i\overline{j}}^{E}=(F_{\alpha\overline{\beta}}^{E}e_{i},e_{j}),F_{\alpha\overline{\beta}}^{L}=(F_{\alpha\overline{\beta}}^{L}\sigma,\sigma),$
and $F^{E},F^{L}$ are the curvature tensors of the Hermitian bundles
$E$ and $L.$

\noindent Since

\noindent \begin{center}
$E\geq_{m}0,m\geq\textrm{min}\{n-k+1,r\}$ and $v_{\overline{\alpha K}i}=0$
for $\alpha\in K,$
\par\end{center}

\noindent we know by definition 3.1

\noindent \begin{center}
$F_{\alpha\overline{\beta}i\overline{j}}^{E}v_{\overline{\alpha K}i}\overline{v_{\overline{\beta K}j}}\geq0.$
\par\end{center}

\noindent From the condition

\noindent \begin{center}
$\sqrt{-1}(\varsigma c(L)-\partial\overline{\partial}\varsigma-\tau^{-1}\partial\varsigma\wedge\overline{\partial}\varsigma)\geq\sqrt{-1}q(\varsigma+\delta)\partial\overline{\partial}\varphi$,
\par\end{center}

\noindent it follows that

$(\left[\sqrt{-1}F^{L\otimes E},\Lambda_{\omega}\right]v,v)\geq(q(\varsigma+\delta)\partial_{\alpha}\partial_{\overline{\beta}}\varphi+\partial_{\alpha}\partial_{\overline{\beta}}\varsigma+\tau^{-1}\partial_{\alpha}\varsigma\overline{\partial_{\beta}\varsigma})v_{\overline{\alpha K}i}\overline{v_{\overline{\beta K}i}}.$

\noindent Substituting the above estimate into (22), we get

$\left\Vert \sqrt{\varsigma}\overline{\partial}^{*}v\right\Vert ^{2}+\left\Vert \sqrt{\varsigma}\overline{\partial}v\right\Vert ^{2}\geq\int_{\Omega}(q(\varsigma+\delta)\partial_{\alpha}\partial_{\overline{\beta}}\varphi+\tau^{-1}\partial_{\alpha}\varsigma\overline{\partial_{\beta}\varsigma})v_{\overline{\alpha K}i}\overline{v_{\overline{\beta K}i}}$

$\ \ \ \ \ \ \ \ \ \ \ \ \ \ \ \ \ \ \ \ \ \ \ \ \ \ \ \ \ \ \ \ \ +2\textrm{Re}(\textrm{grad}^{0,1}\varsigma\lrcorner v,\overline{\partial}^{*}v)dV_{\omega}$

$\ \ \ \ \ \ \ \ \ \ \ \ \ \ \ \ \ \ \ \ \ \ \ \ \ \ \ \ \ \overset{(11)}{=}\int_{\Omega}q(\varsigma+\delta)|\Phi|^{-2}\left\Vert B_{\Phi}A_{v}\right\Vert ^{2}-q(\varsigma+\delta)|\Phi|^{-2}$

$\ \ \ \ \ \ \ \ \ \ \ \ \ \ \ \ \ \ \ \ \ \ \ \ \ \ \ \ \ \ \ \ \ \cdot(F_{A_{v}(\frac{\partial}{\partial z}\otimes\frac{\partial}{\partial\overline{z}_{K}}\otimes\sigma^{*}\otimes e_{i}^{*})\overline{A_{v}(\frac{\partial}{\partial z}\otimes\frac{\partial}{\partial\overline{z}_{K}}\otimes\sigma^{*}\otimes e_{i}^{*})}}^{\textrm{Hom}(E,E^{'})}\Phi,\Phi)$

$\ \ \ \ \ \ \ \ \ \ \ \ \ \ \ \ \ \ \ \ \ \ \ \ \ \ \ \ \ \ \ \ \ \ \ \ \ +\tau^{-1}\partial_{\alpha}\varsigma\overline{\partial_{\beta}\varsigma}v_{\overline{\alpha K}i}\overline{v_{\overline{\beta K}i}}+2\textrm{Re}(\textrm{grad}^{0,1}\varsigma\lrcorner v,\overline{\partial}^{*}v)dV_{\omega}$

$\ \ \ \ \ \ \ \ \ \ \ \ \ \ \ \ \ \ \ \ \ \ \ \ \ \ \ \ \ \geq\left\Vert \sqrt{q(\varsigma+\delta)}B_{\Phi}A_{v}\right\Vert _{|\Phi|^{-2}}^{2}+2\textrm{Re}(\textrm{grad}^{0,1}\varsigma\lrcorner v,\overline{\partial}^{*}v)$

\begin{equation}
+\left\Vert \textrm{grad}^{0,1}\varsigma\lrcorner v\right\Vert _{\tau^{-1}}^{2}.\ \ \end{equation}

Step 3. To deal with the third term in (21), we first do integration
by parts and then apply the H\"{o}lder inequality with an
appropriate parameter.

\noindent Integration by parts yields that

$(e^{-\varphi}\varsigma\Phi^{*}u,\overline{\partial}^{*}v)=(\overline{\partial}(e^{-\varphi}\varsigma\Phi^{*}u),v)$

$\ \ \ \ \ \ \ \ \ \ \ \ \ \ \ \ \ \ \ \ \,=(e^{-\varphi}(-\varsigma\overline{\partial}\varphi\wedge\Phi^{*}u+\varsigma\overline{\partial}\Phi^{*}\wedge u+\overline{\partial}\varsigma\wedge\Phi^{*}u),v)$

$\ \ \ \ \ \ \ \ \ \ \ \ \ \ \ \ \ \ \ \ \,=-(\overline{\partial}\varphi\wedge\Phi^{*}u,v)_{|\Phi|^{-2}\varsigma}+(\overline{\partial}\Phi^{*}\wedge u,v)_{|\Phi|^{-2}\varsigma}$

$\ \ \ \ \ \ \ \ \ \ \ \ \ \ \ \ \ \ \ \ \ \ \ \,+(e^{-\varphi}\Phi^{*}u,\textrm{grad}^{0,1}\varsigma\lrcorner v)$

$\ \ \ \ \ \ \ \ \ \ \ \ \ \ \ \ \ \ \ \ \,=-(\Phi^{*}u,\textrm{grad}^{0,1}\varphi\lrcorner v)_{|\Phi|^{-2}\varsigma}+(\overline{\partial}\Phi^{*}\wedge u,v)_{|\Phi|^{-2}\varsigma}$

\begin{equation}
+(e^{-\varphi}\Phi^{*}u,\textrm{grad}^{0,1}\varsigma\lrcorner v).\ \ \ \ \ \ \ \ \ \ \ \ \end{equation}

\begin{spacing}{0.8}
\noindent In the second equality, we used the condition $\overline{\partial}u=0.$
Now by substituting (15) in lemma 2.3 into (24), we have

\noindent \begin{equation}
(e^{-\varphi}\varsigma\Phi^{*}u,\overline{\partial}^{*}v)=(e^{-\varphi}\Phi^{*}u,\textrm{grad}^{0,1}\varsigma\lrcorner v)+(u,\textrm{Tr}B_{\Phi}A_{v})_{|\Phi|^{-2}\varsigma}.\end{equation}
Lemma 2.2 applied to the homomorphism
\end{spacing}

\begin{center}
$D=B_{\Phi}A_{v}\in\textrm{Hom}(E^{*},\wedge^{n,k}M\otimes L\otimes E^{'}\otimes E^{*})$
\par\end{center}

\noindent gives the following pointwise estimate.

$2|(u,\textrm{Tr}B_{\Phi}A_{v})_{|\Phi|^{-2}\varsigma}|\leq\left\Vert u\right\Vert _{\frac{\varsigma^{2}}{(\varsigma+\delta)|\Phi|^{2}}}^{2}$+$\left\Vert \sqrt{\varsigma+\delta}\textrm{Tr}B_{\Phi}A_{v}\right\Vert _{|\Phi|^{-2}}^{2}$

\begin{equation}
\ \ \ \ \ \ \ \ \ \ \ \ \ \ \ \ \ \ \leq\left\Vert u\right\Vert _{\frac{\varsigma^{2}}{(\varsigma+\delta)|\Phi|^{2}}}^{2}+\left\Vert \sqrt{q(\varsigma+\delta)}B_{\Phi}A_{v}\right\Vert _{|\Phi|^{-2}}^{2}.\end{equation}

\noindent Since

\noindent \begin{center}
$2\textrm{Re}(\textrm{grad}^{0,1}\varsigma\lrcorner v,|\Phi|^{-2}\Phi^{*}u+\overline{\partial}^{*}v)\leq\left\Vert |\Phi|^{-2}\Phi^{*}u+\overline{\partial}^{*}v\right\Vert _{\tau}^{2}+\left\Vert \textrm{grad}^{0,1}\varsigma\lrcorner v\right\Vert _{\tau^{-1}}^{2},$
\par\end{center}

\noindent from (21) (23) (25) (26), it follows that

$\left\Vert |\Phi|^{-2}\Phi^{*}u+\overline{\partial}^{*}v\right\Vert _{\Omega,\varsigma+\tau}^{2}+\left\Vert \overline{\partial}v\right\Vert _{\Omega,\varsigma}^{2}\geq\left\Vert u\right\Vert _{|\Phi|^{-2}\lambda\varsigma}^{2}-\left\Vert u\right\Vert _{\frac{\varsigma^{2}}{(\varsigma+\delta)|\Phi|^{2}}}^{2}$\begin{equation}
\ \ \ \ \ \ \ \ \ \ \ \ \ \ \ \ \ \ \ \ \ \ =\left\Vert u\right\Vert _{\frac{\varsigma(\lambda\delta+\lambda\varsigma-\varsigma)}{(\varsigma+\delta)|\Phi|^{2}}}^{2}.\end{equation}

\noindent This finishes the proof of the main estimate.

\begin{flushright}
$\:$$\square$
\par\end{flushright}

\section*{{\normalsize 4. A Division Theorem for Exact Sequences of Holomorphic
Vector Bundles}}

$\ \ \ $In this section we give a sufficient integrability condition
for the exactness at the level of global holomorphic sections for
exact sequences of holomorphic vector bundles.

We consider a complex of holomorphic vector bundles over $M$, \begin{equation}
E\overset{\Phi}{\rightarrow}E^{'}\overset{\Psi}{\rightarrow}E^{''}\end{equation}
 i.e. $\Phi\in\Gamma(M,\textrm{Hom}(E,E^{'})),\Psi\in\Gamma(M,\textrm{Hom}(E^{'},E^{''}))$
such that $\Psi\circ\Phi=0.$ $E,E^{'},E^{''}$ are assumed to be
endowed with Hermitian structures.

We define for any $x\in M$

\begin{equation}
\mathfrak{\mathcal{E}}(x)=\textrm{min}\{((\Psi^{*}\Psi+\Phi\Phi^{*})\xi,\xi)|\xi\in E_{x}^{'},|\xi|=1\}\end{equation}

\noindent where $\Phi^{*},\Psi^{*}$ are the adjoint of $\Phi$ and
$\Psi$ respectively w.r.t. the given Hermitian structures. By definition,
$0\leq\mathfrak{\mathcal{E}}\in C(M)$ is the smallest eigenvalue
of $\Psi^{*}\Psi+\Phi\Phi^{*}.$ Suppose the complex (28) is exact
at $x\in M$, let $\xi\in E_{x}^{'}$ such that $(\Psi^{*}\Psi+\Phi\Phi^{*})\xi=0,$
then by pairing with $\xi$ we get $\Phi^{*}\xi=0,\Psi\xi=0,$ i.e.
$\xi\in\textrm{Ker}\Phi^{*}\cap\textrm{Ker}\Psi$$=\textrm{Im}\Phi^{\bot}\cap\textrm{Im}\Phi$
which implies $\xi=0.$ Conversely, we assume $\Psi^{*}\Psi+\Phi\Phi^{*}$
is an isomorphism on $E_{x}^{'}$ for some $x\in M$. Since $\textrm{Ker}\Psi$
is invariant under $\Psi^{*}\Psi+\Phi\Phi^{*},$ $\Psi^{*}\Psi+\Phi\Phi^{*}$
also induces an isomorphism on $\textrm{Ker}\Psi.$ Let $\xi\in\textrm{Ker}\Psi,$
there exists some $\eta\in\textrm{Ker}\Psi$ such that $\xi=(\Psi^{*}\Psi+\Phi\Phi^{*})\eta=\Phi\Phi^{*}\eta\in\textrm{Im}\Phi.$
Now we obtain the following useful fact about the function $\mathfrak{\mathcal{E}}$:

\begin{center}
The complex (28) is exact at $x\in M$ if and only if $\mathfrak{\mathcal{E}}(x)>0.$
\par\end{center}

\noindent When the complex (28) is exact, $\Phi^{*}(\Psi^{*}\Psi+\Phi\Phi^{*})^{-1}|_{\textrm{Ker}\Psi}$
is a smooth lifting of $\Phi.$ So it is possible to establish division
theorems by solving a coupled system consisting of

\noindent \begin{center}
$\overline{\partial}g=\overline{\partial}[\Phi^{*}(\Psi^{*}\Psi+\Phi\Phi^{*})^{-1}f]$
and $\Phi g=0$
\par\end{center}

\noindent where $f\in\Gamma(E^{'})$ satisfying $\Psi f=0.$ If $g$
is a solution of this system, then $h\overset{def}{=}\Phi^{*}(\Psi^{*}\Psi+\Phi\Phi^{*})^{-1}f-g\in\Gamma(E)$
and $\Phi h=f.$ If $\Phi$ is surjective and $E^{'}$ is equipped
with the quotient Hermitian structure then $\Psi=0$, $\Phi\Phi^{*}=Id_{E^{'}},$
and the above system reduces to

\noindent \begin{center}
$\overline{\partial}g=\overline{\partial}(\Phi^{*}f)$
\par\end{center}

\noindent on the subbundle $\textrm{Ker}\Phi.$ This key observation
played an important role in both {[}S78{]} and {[}D82{]}. The
difficulty of this method for our case is that $\textrm{Ker}\Phi$ is
no longer a subbundle of $E,$ so it amounts to solving
$\overline{\partial}$-equations for solutions valued in a subsheaf,
it seems that it is not easy to give sufficient conditions for the
solvability of this system.\\

The following lemma reduces our main theorem to the estimate (20).
It was first formulated in {[}S72{]}, the present version is quoted
from {[}V08{]}).\\

\noindent $\mathbf{Lemma4.1.}$ Let $H,H_{0},H_{1},H_{2}$ be Hilbert
spaces, $T:H_{0}\rightarrow H$ be a bounded linear operator, $T_{\ell}:H_{\ell-1}\rightarrow H_{\ell}(\ell=1,2)$
be linear, closed, densely defined operators such that $T_{2}\circ T_{1}=0,$
and let $F\subseteq H$ be a closed subspace such that $T(\textrm{Ker}T_{1})\subseteq F.$
Then for every $f\in F$ the following statements are equivalent

1. there exists at least one $u\in\textrm{Ker}T_{1}$ and $C>0$ such
that $Tu=f$, $\left\Vert u\right\Vert _{H_{0}}\leq C.$

2. $|(g,f)_{H}|^{2}\leq C^{2}(\left\Vert
T^{*}g+T_{1}^{*}v\right\Vert _{H_{0}}^{2}+\left\Vert
T_{2}v\right\Vert _{H_{2}}^{2})$ holds for any $g\in
F,v\in\textrm{Dom(}T_{1}^{*})\cap\textrm{Dom(}T_{2}).$\\

The complex (28) is said to be generically exact if it is exact
outside a subset of measure zero(w.r.t.$dV_{\omega}$) of $M.$\\

\noindent $\mathbf{Theorem4.2.}$ Let $(M,\omega)$ be a K\"{a}hler
manifold and let $E,E^{'},$ $E^{''}$ be Hermitian holomorphic vector
bundles over $M$, $L$ a Hermitian line bundle over $M.$ All the
Hermitian structures may have singularities in a subvariety
$Z\subsetneqq M$ and $\Phi^{-1}(0)\subseteq Z$. Suppose that (28) is
generically exact over $M,$ $M\setminus Z$ is weakly pseudoconvex
and that the following conditions hold on $M\setminus Z$:

1. $E\geq_{m}0,m\geq\textrm{min}\{n-k+1,r\},1\leq k\leq n$;

2. the curvature of $\textrm{Hom}(E,E^{'})$ satisfies

\noindent \begin{center}
$(F_{X\overline{X}}^{\textrm{Hom}(E,E^{'})}\Phi,\Phi)\leq0$ for every
$X\in T^{1,0}M$;
\par\end{center}

3. the curvature of $L$ satisfies

\begin{center}
$\sqrt{-1}(\varsigma c(L)-\partial\overline{\partial}\varsigma-\tau^{-1}\partial\varsigma\wedge\overline{\partial}\varsigma)\geq\sqrt{-1}q(\varsigma+\delta)\partial\overline{\partial}\varphi$.
\par\end{center}

\noindent Then for every $\overline{\partial}$-closed $(n,k-1)$-form
$f$ which is valued in $L\otimes E^{'}$ with $\Psi f=0$ and
$\left\Vert f\right\Vert
_{\frac{\varsigma+\delta}{(\varsigma+\delta)\varsigma\mathfrak{\mathcal{E}}-|\Phi|^{2}\varsigma^{2}}}<+\infty$,
there exists a $\overline{\partial}$-closed $(n,k-1)$-form $h$
valued in $L\otimes E$ such that $\Phi h=f$ and \begin{equation}
\left\Vert h\right\Vert _{\frac{1}{\varsigma+\tau}}\leq\left\Vert
f\right\Vert
_{\frac{\varsigma+\delta}{(\varsigma+\delta)\varsigma\mathfrak{\mathcal{E}}-|\Phi|^{2}\varsigma^{2}}},\end{equation}
where $q=\underset{M\setminus Z}{\textrm{max}}\
\textrm{rank}B_{\Phi},\varphi=\log\left\Vert \Phi\right\Vert $,
$\mathcal{E}$ is the function defined by (29), $0<\varsigma,\tau\in
C^{\infty}(M)$ and $\delta$ is a measurable function on $M$
satisfying
$\mathcal{E}(\varsigma+\delta)\geq||\Phi||^{2}\varsigma.$\\

\noindent Proof. Step 1. Let $\phi\in C^{\infty}(M\setminus Z)$ be
a plurisubharmonic exhaustion function on $M\setminus Z$. For any
$t>0,$ set $\Omega_{t}=\{x\in M\setminus Z|\phi(x)<t\}.$ We know
by definition $\Omega_{t}\Subset M\setminus Z\ \textrm{and}\ \underset{t}{\cup}\Omega_{t}=M\setminus Z$.

\noindent Apart from a subset of $\mathbb{R}$ in $t$ which has measure
zero, $\Omega_{t}$ is a pseudoconvex domain with smooth boundary,
so our main estimate holds on such $\Omega_{t}.$ If we could find
a $\overline{\partial}$-closed section, say $h_{t},$ solving the
equation $\Phi h_{t}=f$ on $\Omega_{t}$ with the estimate

\begin{center}
$\left\Vert h_{t}\right\Vert _{\Omega_{t},\frac{1}{\varsigma+\tau}}\leq\left\Vert f\right\Vert _{\Omega_{t},\frac{\varsigma+\delta}{(\varsigma+\delta)\varsigma\mathfrak{\mathcal{E}}-|\Phi|^{2}\varsigma^{2}}}\leq\left\Vert f\right\Vert _{\frac{\varsigma+\delta}{(\varsigma+\delta)\varsigma\mathfrak{\mathcal{E}}-|\Phi|^{2}\varsigma^{2}}.}$
\par\end{center}

\noindent By setting $h_{t}$ to be zero outside $\Omega_{t}$, we
extend $h_{t}$ to be an element of $L_{n,k-1}^{2}(M,L\otimes E,dV_{\frac{1}{\varsigma+\tau}})$.
The above estimate allows us to apply the compactness argument to
produce on $M$ a $(n,k-1)$-form $h$ valued in $L\otimes E$ as
the weak limit of $\{h_{t}\}$ in $L_{n,k-1}^{2}(M,L\otimes E,dV_{\frac{1}{\varsigma+\tau}})$
which satisfies the equation

\noindent \begin{equation}
\overline{\partial}h=0\:\textrm{outside}\: Z.\end{equation}
and the following inequality

\noindent \begin{center}
$\left\Vert h\right\Vert _{\frac{1}{\varsigma+\tau}}\leq\underset{t\rightarrow+\infty}{\underline{\textrm{lim}}}\left\Vert h_{t}\right\Vert _{\Omega_{t},\frac{1}{\varsigma+\tau}}\leq\left\Vert f\right\Vert _{\frac{\varsigma+\delta}{(\varsigma+\delta)\varsigma\mathfrak{\mathcal{E}}-|\Phi|^{2}\varsigma^{2}}}.$
\par\end{center}

\noindent By using the resulting $L^{2}$-estimate and definition
3.2, we know $h$ has $L_{loc}^{2}$ (w.r.t the Lebesgue measure on
the coordinate chart) coefficients under holomorphic frames over holomorphic
coordinate charts. Since $Z$ is an analytic subset of $M,$ the extension
lemma(D82, lemma 6.9) shows that $\overline{\partial}h=0$ holds on
the whole manifold $M$. Consequently, $\Phi h=f$ on $M$ follows
from solving $\Phi h_{t}=f$ on $\Omega_{t}$ with a $L^{2}$-estimate
for each $t>0.$

$\:$Step 2. In order to make use of lemma 4.1, we introduce the following
Hilbert spaces and densely defined operators.

\begin{center}
$H_{0}=L_{n,k-1}^{2}(\Omega,L\otimes E,dV_{\omega}),\: H=L_{n,k-1}^{2}(\Omega,L\otimes E^{'},dV_{|\Phi|^{-2}}),$
\par\end{center}

\begin{center}
$H_{1}=L_{n,k}^{2}(\Omega,L\otimes E,dV_{\omega}),\ \;\;\;\;\;\;\;\; H_{2}=L_{n,k+1}^{2}(\Omega,L\otimes E,dV_{\omega}),$
\par\end{center}

\begin{center}
$T=\Phi\circ\sqrt{\varsigma+\tau}:H_{0}\rightarrow H,\:\ \;\;\;\;\;\;\; T_{1}=\overline{\partial}\circ\sqrt{\varsigma+\tau}:H_{0}\rightarrow H_{1},$
\par\end{center}

\begin{center}
$T_{2}=\sqrt{\varsigma}\circ\overline{\partial}:H_{1}\rightarrow H_{2},\ \;\;\;\;\;\;\;\;\ \;\;\;\;\;\;\;\;\ \;\;\;\;\;\;\;\;\ \;\;\;\;\ \;\ \ \;\;\;\;\;\;\;\ \;\;\;\;\;\;\;\;\;\;\;\;$
\par\end{center}

\noindent where $\Omega\Subset M\setminus Z$ is a domain satisfying
conditions in lemma 3.3. Then $T:H_{0}\rightarrow H$ is bounded(note
that $\Omega\Subset M\setminus Z$ ), $T_{\ell}:H_{\ell-1}\rightarrow H_{\ell}(\ell=1,2)$
are closed, densely defined and satisfy $T_{2}\circ T_{1}=0.$ It
is easy to see that the adjoint of $T$ is given by

\noindent \begin{center}
$T^{*}u=\sqrt{\varsigma+\tau}|\Phi|^{-2}\Phi^{*}u,u\in H.$
\par\end{center}

\noindent Similarly, from $0<\varsigma,\tau\in C^{\infty}(M)$ we
know

\noindent \begin{center}
$\textrm{Dom}T_{1}^{*}=\textrm{Dom}(\overline{\partial}^{*}),\textrm{Dom}T_{2}=\textrm{Dom}(\overline{\partial})$
\par\end{center}

\noindent and

\noindent \begin{center}
$T_{1}^{*}v=\sqrt{\varsigma+\tau}\overline{\partial}^{*}v,$ $v\in\textrm{Dom}(T_{1}^{*}).$
\par\end{center}

\noindent Define

\noindent \begin{center}
$F=\{u\in H|\Psi u=0,\overline{\partial}u=0\},$
\par\end{center}

\noindent it is easy to see that $T(\textrm{Ker}T_{1})\subseteq F.$
Since $\Psi$ and $\overline{\partial}$ are both closed operators,
\emph{F }is a closed subspace of $H.$

\noindent By the definition of $\mathcal{E}$, we know the following
inequality

\noindent \begin{center}
$|\Phi^{*}u|^{2}=((\Psi^{*}\Psi+\Phi\Phi^{*})u,u)$\begin{equation}
\ \ \ \ \ \ \ \ \ \ \ \ \ \geq\lambda|\Phi|^{2}|u|^{2},\lambda=|\Phi|^{-2}\mathcal{E}.\end{equation}

\par\end{center}

\noindent holds a.e.(w.r.t.$dV_{\omega}$) on $\Omega$ for every
$u\in F.$

\noindent Let $f\in A^{n,k-1}(L\otimes E^{'})$ which is $\overline{\partial}$-closed
and satisfies $\Psi f=0$ then we know by definition $f\in F.$ Since
$\mathcal{E}(\varsigma+\delta)\geq|\Phi|^{2}\varsigma,$ we have

\noindent \begin{center}
$\lambda\delta+\lambda\varsigma\geq\varsigma.$
\par\end{center}

\noindent From the a priori estimate (20) and the density lemma, we
obtain the following inequality

$|(u,f)_{H}|^{2}\leq\left\Vert f\right\Vert _{\frac{\varsigma+\delta}{(\lambda\delta+\lambda\varsigma-\varsigma)\varsigma|\Phi|^{2}}}^{2}\left\Vert u\right\Vert _{\frac{\varsigma(\lambda\delta+\lambda\varsigma-\varsigma)}{(\varsigma+\delta)|\Phi|^{2}}}^{2}$

$\ \ \ \ \ \ \ \ \ \ \ \ \ \leq\left\Vert f\right\Vert _{\frac{\varsigma+\delta}{(\lambda\delta+\lambda\varsigma-\varsigma)\varsigma|\Phi|^{2}}}^{2}(\left\Vert |\Phi|^{-2}\Phi^{*}u+\overline{\partial}^{*}v\right\Vert _{\varsigma+\tau}^{2}+\left\Vert \overline{\partial}v\right\Vert _{\varsigma}^{2})$

$\ \ \ \ \ \ \ \ \ \ \ \ \ =\left\Vert f\right\Vert _{\frac{\varsigma+\delta}{(\lambda\delta+\lambda\varsigma-\varsigma)\varsigma|\Phi|^{2}}}^{2}(\left\Vert \sqrt{\varsigma+\tau}|\Phi|^{-2}\Phi^{*}u+\sqrt{\varsigma+\tau}\overline{\partial}^{*}v\right\Vert ^{2}$

$\ \ \ \ \ \ \ \ \ \ \ \ \ \ \ \ +\left\Vert \sqrt{\varsigma}\overline{\partial}v\right\Vert ^{2})$

$\ \ \ \ \ \ \ \ \ \ \ \ \ =\left\Vert f\right\Vert _{\frac{\varsigma+\delta}{(\lambda\delta+\lambda\varsigma-\varsigma)\varsigma|\Phi|^{2}}}^{2}(\left\Vert T^{*}u+T_{1}^{*}v\right\Vert _{H_{0}}^{2}+\left\Vert T_{2}v\right\Vert _{H_{2}}^{2})$

\begin{spacing}{0.7}
\noindent holds for any $u\in F,v\in\textrm{Dom}(T_{1}^{*})\cap\textrm{Dom(}T_{2}).$
Note that the condition $\lambda\delta+\lambda\varsigma\geq\varsigma$
is needed for the first inequality. Hence we know by lemma 5.1 that
there exist at least one $h^{'}\in\textrm{Ker}T_{1}$ such that
\end{spacing}

\noindent \begin{center}
$Th^{'}=f$ and $\left\Vert h^{'}\right\Vert _{H_{0}}^{2}\leq\left\Vert f\right\Vert _{\frac{\varsigma+\delta}{(\lambda\delta+\lambda\varsigma-\varsigma)\varsigma|\Phi|^{2}}}^{2}.$
\par\end{center}

\noindent Letting $h=\sqrt{\varsigma+\tau}h^{'},$ we have

\noindent \begin{center}
$\Phi h=f$, $\left\Vert h\right\Vert _{\frac{1}{\varsigma+\tau}}^{2}\leq\left\Vert f\right\Vert _{\frac{\varsigma+\delta}{(\lambda\delta+\lambda\varsigma-\varsigma)\varsigma|\Phi|^{2}}}^{2}$.
\par\end{center}

\noindent Replacing $\lambda$ by $|\Phi|^{-2}\mathcal{E}$ completes
the proof.

\begin{spacing}{0.7}
\noindent \begin{flushright}
$\ $$\square$
\par\end{flushright}
\end{spacing}

\noindent $\mathbf{Remarks.}$ (i) If $M$ is weakly pseudoconvex and
$Z=\eta^{-1}(0)\subsetneqq M$ where $\eta$ is a holomorphic function
on $M$, then $M\setminus Z$ is weakly pseudoconvex. Let $\psi\in
C^{\infty}(M)$ be a plurisubharmonic exhaustion function on $M$. It
is easy to see that $\phi:=\psi+|\eta|^{-1}$ is a plurisubharmonic
exhaustion function on $M\setminus Z$. (ii) If $M$ is a Stein
manifold (or more generally, an essentially Stein manifold, see
{[}V08{]}) and $Z$ is an analytic hypersurface, then $M\setminus Z$
is a Stein manifold(or an essentially Stein manifold). (iii) When
$\Phi$ is not identically zero, one can always find an analytic
hypersurface $Z$ such that $\Phi^{-1}(0)\subseteq Z$.

$\,$

If we choose $\varsigma$ to be a positive constant, the third condition
in theorem 4.2 will be independent of the function $\tau.$ By this
observation, we have the following corollary.

$\:$

\noindent $\mathbf{Corollary4.3.}$ If the condition 3 in theorem
4.2 is replaced by

\begin{equation}
\sqrt{-1}c(L)\geq\sqrt{-1}q(|\Phi|^{2}\mathcal{E}^{-1}+1)\partial\overline{\partial}\varphi,\end{equation}

\noindent then for every $\overline{\partial}$-closed $(n,k-1)$-form
$f$ which is valued in $L\otimes E^{'}$ with

\noindent \begin{center}
$\Psi f=0$ and $\left\Vert f\right\Vert _{\frac{\mathcal{E}+|\Phi|^{2}}{\mathcal{E}^{2}}}<+\infty$
\par\end{center}

\noindent there is a $\overline{\partial}$-closed $(n,k-1)$-form
$h$ valued in $L\otimes E$ such that $\Phi h=f$ and the following
estimate holds \begin{equation}
\left\Vert h\right\Vert \leq\left\Vert f\right\Vert _{\frac{\mathcal{E}+|\Phi|^{2}}{\mathcal{E}^{2}}}.\end{equation}

\noindent Proof. Set

\noindent \begin{center}
$\varsigma=1,\tau=\textrm{constant}>0,\delta=|\Phi|^{2}\mathcal{E}^{-1},$
\par\end{center}

\noindent it is easy to see that $\mathcal{E}(\varsigma+\delta)\geq|\Phi|^{2}\varsigma.$

\noindent Hence, we get from theorem 4.2 that for every $\overline{\partial}$-closed
$(n,k-1)$-form $f$ which is valued in $L\otimes E^{'}$and satisfies

\noindent \begin{center}
$\Psi f=0,$ $\left\Vert f\right\Vert _{\frac{\varsigma+\delta}{(\varsigma+\delta)\varsigma\mathfrak{\mathcal{E}}-|\Phi|^{2}\varsigma^{2}}}=\left\Vert f\right\Vert _{\frac{\mathcal{E}+|\Phi|^{2}}{\mathcal{E}^{2}}}<+\infty$,
\par\end{center}

\noindent there is at least a $\overline{\partial}$-closed $(n,k-1)$-form
$h_{\tau}$ valued in $L\otimes E$ such that

\noindent \begin{center}
$\Phi h_{\tau}=f$ and $\left\Vert h_{\tau}\right\Vert _{\frac{1}{1+\tau}}\leq\left\Vert f\right\Vert _{\frac{\varsigma+\delta}{(\varsigma+\delta)\varsigma\mathfrak{\mathcal{E}}-|\Phi|^{2}\varsigma^{2}}}=\left\Vert f\right\Vert _{\frac{\mathcal{E}+|\Phi|^{2}}{\mathcal{E}^{2}}}.$
\par\end{center}

\noindent From the estimate given above, it follows that

\noindent \begin{center}
$\left\Vert h_{\tau}\right\Vert =\sqrt{1+\tau}\left\Vert h_{\tau}\right\Vert _{\frac{1}{1+\tau}}$
$\leq\sqrt{1+\tau}\left\Vert f\right\Vert _{\frac{\mathcal{E}+|\Phi|^{2}}{\mathcal{E}^{2}}}.$
\par\end{center}

The above estimate shows that $\{h_{\tau}\}_{1>\tau>0}$ bounded in
$L_{n,k-1}^{2}(M,L\otimes E,dV_{\omega})$, so we get a weak limit
$h$ of $\{h_{\tau}\}_{\tau>0}$ in $L_{n,k-1}^{2}(M,L\otimes E,dV_{\omega})$
when $\tau\rightarrow0.$ It is easy to see that the resulting section
$h$ is $\overline{\partial}$-closed on $M$ and $\Phi h=f.$ The
$L^{2}$-estimate of $h_{\tau}$ implies that

\noindent \begin{center}
$\left\Vert h\right\Vert \leq\underset{\tau\rightarrow0}{\underline{\textrm{lim}}}\left\Vert h_{\tau}\right\Vert $
\par\end{center}

\noindent \begin{center}
$\ \ \ \ \ \ \ \ \ \ \ \ \ \ \ \ \ \ \ \ \ \leq\underset{\tau\rightarrow0}{\underline{\textrm{lim}}}\sqrt{1+\tau}\left\Vert f\right\Vert _{\frac{\mathcal{E}+|\Phi|^{2}}{\mathcal{E}^{2}}}$
\par\end{center}

\noindent \begin{center}
$\ \ \ \ \ \ =\left\Vert f\right\Vert _{\frac{\mathcal{E}+|\Phi|^{2}}{\mathcal{E}^{2}}}$
\par\end{center}

\begin{spacing}{0.7}
\noindent which completes the proof of corollary 4.3.
\end{spacing}

\begin{spacing}{0.7}
\noindent \begin{flushright}
$\ $$\square$
\par\end{flushright}
\end{spacing}

\noindent $\mathbf{Remarks.}$ (i) The condition $(F_{X\overline{X}}^{\textrm{Hom}(E,E^{'})}\Phi,\Phi)\leq0$
is needed to handle the second term in (11). We recall that the curvature
of the Chern connection of a Hermitian holomorphic vector bundle is
semi-negative in the sense of Griffiths(Nakano) if and only if it
is $1$-tensor($\textrm{min}\{n,r\}$-tensor) semi-negative. Hence
a sufficient condition for $(F_{X\overline{X}}^{\textrm{Hom}(E,E^{'})}\Phi,\Phi)\leq0$
is given by(since we always assume $E\geq_{m}0$ for some positive
integer $m$): $E^{'}$ is semi-negative in the sense of Griffiths.

\noindent (ii) In particular, if the underlying manifold $M$ is
assumed to be strongly pseudoconvex(corollary 4.5 (ii)) then one can
always endow the holomorphic vector bundles over $M,$ $E$ and
$E^{'}$ with Hermitian structures such that $E$ is semi-positive in
the sense of Nakano and $E^{'}$ is semi-negative in the sense of
Griffiths. So our curvature conditions 1 and 2 are satisfied
automatically by such Hermitian structures.

\noindent (iii) For these homomorphisms in the Koszul complex, due
to the identity (50), the condition $(F_{X\overline{X}}^{\textrm{Hom}(E,E^{'})}\Phi,\Phi)\leq0$
holds provided the $E$ is semi-positive in sense of Griffiths. The
Koszul complex provides a series of homomorphisms which are not generically
surjective. The generically surjective case has been extensively investigated
in {[}S78{]} and {[}D82{]}.

$\,$

We can derive from corollary 4.3 the following results.

$\:$

\noindent $\mathbf{Corollary4.4.}$ Results in corollary 4.3 hold
with the condition 2 assumed there replaced by the condition that
$E^{'}$ is semi-negative in the sense of Griffiths.

$\:$

\noindent $\mathbf{Corollary4.5.}$ Besides the conditions in corollary
4.3, we also assume that (28) is exact on the whole manifold $M$.
Then we have

\begin{spacing}{0.8}
\noindent (i) For every $\overline{\partial}$-closed $(n,k-1)$-form
$f$ which is valued in $L\otimes E^{'}$ and locally
square-integrable on $M,$ if $\Psi f=0$ then there exists a
$\overline{\partial}$-closed $h\in L_{n,k-1}^{2}(M,L\otimes
E,dV_{\omega})$ such that $\Phi h=f.$ In particular, if $E$ is
semi-positive in the sense of Nakano, then the induced sequence on
global section
\end{spacing}

\begin{spacing}{0.8}
\noindent \begin{center}
\begin{equation}
\Gamma(M,K_{M}\otimes L\otimes E)\rightarrow\Gamma(M,K_{M}\otimes L\otimes E^{'})\rightarrow\Gamma(M,K_{M}\otimes L\otimes E^{''})\end{equation}

\par\end{center}
\end{spacing}

\noindent is exact where $K_{M}$ is the canonical bundle of $M.$

\begin{spacing}{0.8}
\noindent Moreover, if $\Phi$ is surjective then it induces a surjective
homomorphism on cohomology groups:
\end{spacing}

\begin{spacing}{0.8}
\noindent \begin{center}
\begin{equation}
\Phi:H^{n,k-1}(M,L\otimes E)\rightarrow H^{n,k-1}(M,L\otimes E^{'}).\end{equation}

\par\end{center}
\end{spacing}

\noindent (ii) If $(M,\omega)$ is strongly pseudoconvex and $\Phi\in\Gamma(M,\textrm{Hom}(E,E^{'}))$
is nonvanishing then (i) holds without assuming the curvature conditions
1-3.

$\:$

\noindent Proof. The proof consisting of using appropriate weight
functions to modify the given Hermitian structure on $L$ to control
the $L^{2}$-norm and curvature.

\noindent (i) Let $0<\phi\in \rm{PSH}(M)\cap C^{\infty}(M)$ be an
exhaustion function on $M.$ $ $Set $\Omega_{t}=\{z\in
M|\phi(z)<t\},t\in\mathbb{R}$ then $\Omega_{t}\Subset M$ and
$\underset{t}{\cup}\Omega_{t}=M.$ Since (28) is assumed to be exact
on $M,$ then we have $\mathcal{E}(x)>0$ for every $x\in M.$ Given
$f\in A^{n,k-1}(L\otimes E^{'})$, we can define a positive number
for each $\ell=0,1,2\cdots$

\begin{center}
$\delta_{\ell}=\textrm{sup}\{\frac{\mathcal{E}(x)+|\Phi(x)|^{2}}{\mathcal{E}(x)^{2}}|x\in\Omega_{\ell+1}\setminus\Omega_{\ell}\}\int_{\Omega_{\ell+1}\setminus\Omega_{\ell}}|f|^{2}dV_{\omega}\in[0,+\infty).$
\par\end{center}

\begin{spacing}{0.7}
\noindent We choose an increasing convex function $\eta\in C^{\infty}(\mathbb{R})$
such that
\end{spacing}

\begin{spacing}{0.7}
\noindent \begin{center}
\begin{equation}
\eta(\ell)\geq\textrm{log}(2^{\ell}\delta_{\ell})\ \textrm{for}\ \ell=0,1,2,\cdots\end{equation}
 $ $
\par\end{center}
\end{spacing}

\noindent and set $\psi=\eta\circ\phi.$ Then $\psi\in\rm{PSH}(M)\cap
C^{\infty}(M)$ and $e^{-\psi}h_{L}$ defines a singular Hermitian
structure on the line bundle $L$ where we denote by $h_{L}$ the
given Hermitian structure on $L.$

\noindent It is easy to see that on $M\setminus\Phi^{-1}(0)$ the
curvature of $e^{-\psi}h_{L}$ satisfies

\noindent \begin{center}
$\sqrt{-1}c(L,e^{-\psi}h_{L})=\sqrt{-1}(\partial\overline{\partial}\psi+c(L,h_{L})\geq\sqrt{-1}q(|\Phi|^{2}\mathcal{E}^{-1}+1)\partial\overline{\partial}\varphi.$
\par\end{center}

\noindent By the construction of $\psi,$ we get the following estimate
of the $L^{2}$-norm of $f$ where the left hand side is computed
by using the new Hermitian structure $e^{-\psi}h_{L}$ on $L.$

\begin{center}
$\left\Vert f\right\Vert _{\frac{\mathcal{E}+|\Phi|^{2}}{\mathcal{E}^{2}}}^{2}=\int_{M}|f|^{2}\frac{\mathcal{E}+|\Phi|^{2}}{\mathcal{E}^{2}}e^{-\psi}dV_{\omega}$
\par\end{center}

$\ \ \ \ \ \ \ \ \ \ \ \ \ \ \ \ \ \ \ \ \ \ \ \ \ \ \ \ \ \ \ \ \ \ \ \ ={\displaystyle \sum_{\ell\geq0}}\int_{\Omega_{\ell+1}\setminus\Omega_{\ell}}|f|^{2}\frac{\mathcal{E}+|\Phi|^{2}}{\mathcal{E}^{2}}e^{-\psi}dV_{\omega}$\begin{equation}
\ \ \ \ \ \ \ \ \ \ \ \overset{(37)}{\leq}\sum_{\ell\geq0}2^{-\ell}=2<+\infty.\end{equation}

\noindent From (38) and corollary 4.3, we get a $\overline{\partial}$-closed
section $h\in L_{n,k-1}^{2}(M,L\otimes E,dV_{\omega})$ such that
$\Phi h=f$ provided $\overline{\partial}f=0.$ Consequently, by using
the De Rham-Weil isomorphism theorem, we know that if $\Psi=0$ then
the induced homomorphism $\Phi:H^{n,k-1}(M,L\otimes E)\rightarrow H^{n,k-1}(M,L\otimes E^{'})$
is surjective.

\noindent In the case of $k=1$, our condition 1 becomes that $E$
is nonnegative in the sense of Nakano. From the ellipticity of $\overline{\partial}$
( also due to the condition that $k=1$), it follows that the induced
sequence (35) is still exact.

(ii) As $M$ is strongly pseudoconvex, one can modify the given
Hermitian structures for $E$ and $E^{'}$ such that $E$ is
semi-positive in the of Nakano and $E^{'}$ is semi-negative in the
sense of Griffiths. With such Hermitian structures, we get the
desired curvature conditions 1 and 2. Next, we multiply the given
Hermitian structure on $L$ by certain weight to make the new
Hermitian structure satisfy condition 3. Let $\phi\in C^{\infty}(M)$
be a strictly plurisubharmonic exhaustion function of $M.$

Set

\noindent \begin{center}
$\lambda(x):=$the smallest eigenvalue of $\sqrt{-1}\partial\overline{\partial}\phi(x)$
w.r.t. the metric $\omega$,
\par\end{center}

\noindent \begin{center}
$\ \ \ \ \ \ \mu(x):=$the smallest eigenvalue of $\sqrt{-1}c(L,h_{L})(x)$
w.r.t. the metric $\omega$,
\par\end{center}

\noindent \begin{center}
$\Lambda(x):=$ the largest eigenvalue of $\sqrt{-1}\partial\overline{\partial}\varphi(x)$
w.r.t. the metric $\omega$,
\par\end{center}

\noindent then it is easy to see that $\lambda,\mu,\Lambda\in C(M)$
and

\noindent \begin{center}
$\sqrt{-1}\partial\overline{\partial}\phi\geq\lambda\omega,\sqrt{-1}c(L,h_{L})\geq\mu\omega,\sqrt{-1}\partial\overline{\partial}\varphi\leq\Lambda\omega.$
\par\end{center}

Since $\phi$ is strictly plurisubharmonic, we know $\lambda>0$ on
$M.$ We can therefore define a function $\sigma:\mathbb{\mathbb{R}\rightarrow R}$
as follows

\noindent \begin{equation}
\sigma(t)=\textrm{sup}\{\frac{q(|\Phi|^{2}(x)\mathcal{E}(x)^{-1}+1)\Lambda(x)-\mu(x)}{\lambda(x)}|x\in\Omega_{t}\},t\in\mathbb{R}.\end{equation}

\noindent For this $\sigma(t),$ one can always find a $\chi\in C^{\infty}[0,+\infty)$
such that \begin{equation}
\chi^{'}(t)\geq\textrm{max}\{\sigma(t),0\},\chi^{''}(t)\geq0,t>0.\end{equation}

Now endow the line bundle $L$ with a new Hermitian structure $e^{-\psi_{1}}h_{L}$
where $\psi_{1}=\chi\circ\phi.$ It is obvious that

$\,\;\sqrt{-1}c(L,e^{-\psi_{1}}h_{L})=\sqrt{-1}(\partial\overline{\partial}\psi_{1}+c(L,h_{L}))$

$\ \ \ \ \ \ \ \ \ \ \ \ \ \ \ \ \ \ \ \ \ \ \ \ \ =\sqrt{-1}(\chi^{'}\circ\phi\partial\overline{\partial}\phi+\chi^{''}\circ\phi\partial\phi\wedge\overline{\partial}\phi+c(L,h_{L}))$

$\ \ \ \ \ \ \ \ \ \ \ \ \ \ \ \ \ \ \ \ \ \ \ \ \overset{(40)}{\geq}\chi^{'}\circ\phi\lambda\omega+\sqrt{-1}c(L,h_{L})$

$\ \ \ \ \ \ \ \ \ \ \ \ \ \ \ \ \ \ \ \ \ \ \ \ \overset{\lambda>0}{\geq}\sigma\circ\phi\lambda\omega+\sqrt{-1}c(L,h_{L})$

$\ \ \ \ \ \ \ \ \ \ \ \ \ \ \ \ \ \ \ \ \ \ \ \ \overset{(39)}{\geq}(q(|\Phi|^{2}\mathcal{E}^{-1}+1)\Lambda-\mu)\omega+\sqrt{-1}c(L,h_{L})$\begin{equation}
\geq\sqrt{-1}q(|\Phi|^{2}\mathcal{E}^{-1}+1)\partial\overline{\partial}\varphi.\ \ \ \ \ \end{equation}

\noindent Conclusion (ii) follows from (41) and conclusion (i), this
finishes the proof of corollary 4.5.

\begin{spacing}{0.7}
\begin{flushright}
$\square$
\par\end{flushright}
\end{spacing}

Given a holomorphic section $\Phi\in\Gamma(M,\textrm{Hom}(E,E^{'})),$
there is the following exact sequence of sheaves over $M$

\begin{center}
$\mathcal{O}(E)\rightarrow\mathcal{O}(E^{'})\rightarrow\mathcal{O}(E^{'})/\textrm{Im}\Phi\rightarrow0$
\par\end{center}

\noindent where ${\rm Im}\Phi$ is the image of the induced
sheaf-homomorphism
$\Phi:\mathcal{O}(E)\rightarrow\mathcal{O}(E^{'}).$ Generally, the
quotient sheaf $\mathcal{O}(E^{'})/\textrm{Im}\Phi$ is never locally
free, so this case does not fit into the general setting established
by theorem 4.2, corollary 4.3 and corollary 4.4.

However we can modify definition of the function $\mathfrak{\mathcal{E}}$
to make the same argument works for this situation. To this end, we
have to introduce the following function $\mathfrak{\mathcal{E}}_{1}$
(instead of $\mathfrak{\mathcal{E}}$) by using a single homomorphism
$\Phi.$ We define for any $x\in M,$

\begin{flushright}
$\mathfrak{\mathcal{E}}_{1}(x)=\textrm{min}\{((\Psi^{*}\Psi+\Phi\Phi^{*})\xi,\xi)|\xi\in E_{x}^{'},|\xi|=1\}$
$\ \ \ \ \ \ \ \ \ \ \ \ \ \ (29)^{'}$
\par\end{flushright}

\noindent where $\Psi$ is the orthogonal projection from $E_{x}^{'}$
onto the subspace $\Phi(E_{x})^{\bot}.$

\noindent It is obvious that $\mathfrak{\mathcal{E}}_{1}$ is positive
everywhere. From the fact that $|\Psi\xi|=\textrm{inf}\{|\xi+\Phi(\eta)||\eta\in E_{x}\}$
for every $\xi\in E_{x}^{'}$ and our definition $(29)^{'}$ we know
the function $\mathfrak{\mathcal{E}_{1}}$ is upper semi-continuous
and therefore measurable.

$\,$

Similar to theorem 4.2, we have the following result about the division
problem for a single holomorphic homomorphism.

$\,$

\noindent $\mathbf{Theorem4.2^{'}.}$ Let $(M,\omega)$ be a
K\"{a}hler manifold and let $E,E^{'},$ $E^{''}$ be Hermitian
holomorphic vector bundles over $M$, $L$ a Hermitian line bundle
over $M.$ All the Hermitian structures may have singularities in a
subvariety $Z\subsetneqq M$ and $\Phi^{-1}(0)\subseteq Z$. Suppose
that $M\setminus Z$ is weakly pseudoconvex and that the following
conditions hold on $M\setminus Z$:

1. $E\geq_{m}0,m\geq\textrm{min}\{n-k+1,r\},1\leq k\leq n$;

2. the curvature of $\textrm{Hom}(E,E^{'})$ satisfies

\noindent \begin{center}
$(F_{X\overline{X}}^{\textrm{Hom}(E,E^{'})}\Phi,\Phi)\leq0$ for every
$X\in T^{1,0}M$;
\par\end{center}

3. the curvature of $L$ satisfies

\begin{center}
$\sqrt{-1}(\varsigma c(L)-\partial\overline{\partial}\varsigma-\tau^{-1}\partial\varsigma\wedge\overline{\partial}\varsigma)\geq\sqrt{-1}q(\varsigma+\delta)\partial\overline{\partial}\varphi$.
\par\end{center}

\noindent Then for every $\overline{\partial}$-closed $(n,k-1)$-form
$f$ which is valued in $L\otimes E^{'}$ and satisfies

\noindent \begin{center}
$f(x)\in\Phi(E_{x})$ for a.e. $x\in M$ and $\left\Vert f\right\Vert _{\frac{\varsigma+\delta}{(\varsigma+\delta)\varsigma\mathfrak{\mathfrak{\mathcal{E}}_{1}}-|\Phi|^{2}\varsigma^{2}}}<+\infty$,
\par\end{center}

\noindent there exists a $\overline{\partial}$-closed $(n,k-1)$-form
$h$ valued in $L\otimes E$ such that $\Phi h=f$ and $\left\Vert h\right\Vert _{\frac{1}{\varsigma+\tau}}\leq\left\Vert f\right\Vert _{\frac{\varsigma+\delta}{(\varsigma+\delta)\varsigma\mathfrak{\mathfrak{\mathcal{E}}_{1}}-|\Phi|^{2}\varsigma^{2}}},$
where $q=\underset{M\setminus Z}{\textrm{max}}\:\textrm{rank}B_{\Phi},\varphi=\log\left\Vert \Phi\right\Vert ^{2}$,
$\mathfrak{\mathcal{E}}_{1}$ is the function defined by $(29)^{'}$,
$0<\varsigma,\tau\in C^{\infty}(M)$ and $\delta$ is a measurable
function on $M$ satisfying $\mathfrak{\mathcal{E}}_{1}(\varsigma+\delta)\geq||\Phi||^{2}\varsigma.$

$\ $

\noindent $\mathbf{Remark.}$ The results parallel to corollaries
4.3-4.4 can be easily derived from theorem $4.2^{'}$. Since the function
$\mathfrak{\mathcal{E}}_{1}$ defined by (29)$^{'}$ is only upper
semi-continuous, it can't be locally bounded from below by positive
constants. So we don't have the result parallel to corollary 4.5.

\section*{{\normalsize 5. Applications to Koszul Complex}}

$\ \ \ $In this section we apply results obtained in section 4 to
the special case of generalized Koszul complex.

Let $M$ be a complex manifold and $E$ be a holomorphic vector bundle
of rank $r$ over $M.$ The Koszul complex associated to a section
$s\in\Gamma(E^{*})$ is defined as follows

\begin{equation}
\textrm{det}E\overset{d_{r}}{\rightarrow}\wedge^{r-1}E\overset{d_{r-1}}{\rightarrow}\cdots\overset{d_{1}}{\rightarrow}\mathcal{O}_{M}\overset{d_{0}}{\rightarrow}0\end{equation}

\noindent where the boundary operators are given by the interior product
\begin{equation}
d_{p}=s\lrcorner,1\leq p\leq r.\end{equation}
 It forms a complex since we have $d_{p-1}\circ d_{p}=0$ for $1\leq p\leq r.$

In particular, if we set $E=T^{*1,0}M$ and $s=X,$ a holomorphic
vector field on $M,$ then the complex (42) is given by

\begin{center}
$K_{M}\overset{X\lrcorner}{\rightarrow}\wedge^{n-1}T^{*1,0}M\overset{X\lrcorner}{\rightarrow}\cdots\rightarrow T^{*1,0}M\overset{X\lrcorner}{\rightarrow}\mathcal{O}_{M}\overset{X\lrcorner}{\rightarrow}0$
\par\end{center}

\noindent which recovers the usual notion of the Koszul complex associated
to a vector field $X$ on a complex manifold $M.$

As before, in order to handle the curvature term in the Bochner formula
we consider the following complex associated to (42):

\begin{center}
$L\otimes\textrm{det}E\overset{d_{r}}{\rightarrow}L\otimes\wedge^{r-1}E\overset{d_{r-1}}{\rightarrow}\cdots\overset{d_{1}}{\rightarrow}L\overset{d_{0}}{\rightarrow}0,$
\par\end{center}

\noindent where $L$ is a holomorphic line bundle over $M.$

We start with improving the estimate in lemma 3.3 for

\begin{center}
$\Phi=s\lrcorner\in\Gamma(M,\textrm{Hom}(\wedge^{p}E,\wedge^{p-1}E)$
\par\end{center}

\noindent $1\leq p\leq r.$ In the following discussion, $E$ and
$L$ are endowed with Hermitian structures. Let \emph{$\{e_{1},\cdots,e_{r}\}$
}be a local holomorphic frame of\emph{ $E$} and let $\{e_{_{1}}^{*},\cdots,e_{r}^{*}\}$
be its dual frame. We conclude form the definition

\begin{center}
$\Phi^{*}=\theta\wedge$
\par\end{center}

\begin{flushleft}
where \begin{equation}
\theta=\overline{g}_{i}h^{\overline{i}j}e_{j}\:\textrm{and}\: s=g_{i}e_{i}^{*}.\end{equation}

\par\end{flushleft}

\noindent By choosing a local frame \emph{$\{e_{1},\cdots,e_{r}\}$
}normal at a given point $x\in M\setminus s^{-1}(0)$ such that

\begin{center}
$e_{1}^{*}(x)=\frac{s}{|s|}(x),$
\par\end{center}

\begin{flushleft}
then we get
\par\end{flushleft}

$||\Phi||^{2}(x)=\underset{i_{1}<\cdots<i_{p}}{\sum}|\Phi(e_{i_{1}}\wedge\cdots\wedge e_{i_{p}})|^{2}$

\begin{equation}
=\underset{1<i_{2}<\cdots<i_{p}}{\sum}||s|e_{i_{2}}\wedge\cdots\wedge e_{i_{p}}|^{2}=\binom{r}{p-1}|s|^{2}(x).\ \ \end{equation}
Since for any $\xi\in\wedge^{p-1}E_{x},x\in M,$ we have

\begin{center}
$\theta\wedge s\lrcorner\xi+s\lrcorner\theta\wedge\xi=|s|^{2}\xi,$
\par\end{center}

\noindent so the function $\mathcal{E}$(in (29)) is given by

\noindent \begin{center}
$\mathcal{E}(x)=|s(x)|^{2}.$
\par\end{center}

\noindent This implies that the complex (42) is exact at $x\in M$
if and only if $s(x)\neq0.$

We will denote by $B_{s}$ the second fundamental form of the line
bundle in $E^{*}$ generated by $s$ over $M\setminus s^{-1}(0),$
i.e.\begin{equation}
B_{s}(X)=(\nabla_{X}^{E^{*}}s)^{\perp}\end{equation}

\noindent where $\nabla^{E^{*}}$ is the Chern connection on $E^{*},$
$X\in T_{x}^{1,0}M,x\in M\setminus s^{-1}(0).$

\noindent $\ $

\noindent $\mathbf{Lemma5.1.}$ We have the following relations between
$s$ and the associated homomorphism $\Phi=s\lrcorner$:

\begin{equation}
B_{\Phi}=B_{s}\lrcorner\end{equation}

\begin{equation}
\left\Vert B_{\Phi}A\right\Vert ^{2}=\binom{r}{p-1}\left\Vert B_{s}A\right\Vert ^{2}\end{equation}

\begin{equation}
\textrm{Tr}B_{\Phi}A=\textrm{Tr}B_{s}A\end{equation}

\begin{equation}
(F_{X\overline{X}}^{\textrm{Hom}(\wedge^{p}E,\wedge^{p-1}E)}\Phi,\Phi)=\binom{r}{p-1}(F_{X\overline{X}}^{E^{*}}s,s)\end{equation}

\noindent where $X\in T_{x}^{1,0}M,x\in M\setminus s^{-1}(0),A\in\textrm{Hom}(\wedge^{n,k-1}TM\otimes L^{*}\otimes E^{*},T^{1,0}M)$,
$\textrm{Tr}B_{s}A$ is defined by (12) with $\rho$ being the interior
product and $F^{E^{*}}$ is the curvature of the induced Chern connection
on $E^{*},1\leq k\leq n.$

$\ $

\noindent Proof. We first choose holomorphic coordinates and frames
$\left\{ z_{1},\cdots,z_{n}\right\} ,$ \emph{$\left\{ e_{1},\cdots,e_{r}\right\} ,$}
\{$\sigma\}$ which are normal at a given point $x\in M\setminus s^{-1}(0).$

\noindent For every $\xi\in\wedge^{p}E_{x}$ we have by the definition
of $B$

\noindent $B(X)\cdot\xi=\nabla_{X}^{\textrm{Hom}(\wedge^{p}E,\wedge^{p-1}E)}(\Phi)\cdot\xi$

\noindent $\ \ \ \ \ \ \ \ \ \ \ \ =\nabla_{X}^{\wedge^{p-1}E}(\Phi\xi)-\Phi(\nabla_{X}^{\wedge^{p}E}\xi)$

\noindent $\ \ \ \ \ \ \ \ \ \ \ \ =\nabla_{X}^{\wedge^{p-1}E}(s\lrcorner\xi)-s\lrcorner(\nabla_{X}^{\wedge^{p}E}\xi)=(\nabla_{X}^{E^{*}}s)\lrcorner\xi.$

\noindent Combining this equality with (45) and using the assumption
that $\left\{ e_{1},\cdots,e_{r}\right\} $ is normal at $x,$ we
obtain

\noindent $(B(X),\Phi)=(B(X)\cdot e_{i_{1}}\wedge\cdots\wedge e_{i_{p}},\Phi(e_{i_{1}}\wedge\cdots\wedge e_{i_{p}})$

\noindent $\ \ \ \ \ \ \ \ \ \ \ \ \ \ =(\nabla_{X}^{E^{*}}s\lrcorner(e_{i_{1}}\wedge\cdots\wedge e_{i_{p}}),$$s\lrcorner(e_{i_{1}}\wedge\cdots\wedge e_{i_{p}}))$

\noindent $\ \ \ \ \ \ \ \ \ \ \ \ \ \ =(\nabla_{X}^{\wedge^{p-1}E}(s\lrcorner(e_{i_{1}}\wedge\cdots\wedge e_{i_{p}})),$$s\lrcorner(e_{i_{1}}\wedge\cdots\wedge e_{i_{p}}))$

\noindent $\ \ \ \ \ \ \ \ \ \ \ \ \ \ =X(|s\lrcorner(e_{i_{1}}\wedge\cdots\wedge e_{i_{p}})|^{2})$

\noindent $\ \ \ \ \ \ \ \ \ \ \ \ \ \ =X(|\Phi|^{2})=\binom{r}{p-1}X(|s|^{2}).$

\noindent Now it follows from definition (8) that

\noindent $B_{\Phi}(X)\cdot\xi=B(X)\cdot\xi-(B(X),\Phi)\frac{\Phi(\xi)}{|\Phi|^{2}}$

\noindent $\ \ \ \ \ \ \ \ \ \ \ \ \ \ =(\nabla_{X}^{E^{*}}s-\frac{X(|s|^{2})}{|s|^{2}}s)\lrcorner\xi.$

\noindent $\ \ \ \ \ \ \ \ \ \ \ \ \ \ =(\nabla_{X}^{E^{*}}s-\frac{(\nabla_{X}^{*}s,s))}{|s|^{2}}s)\lrcorner\xi=B_{s}(X)\lrcorner\xi.$

\noindent For increasing multi-indices $K,I$ with $|K|=k-1,|I|=p,$
we denote $X_{\overline{K}I}=A(\frac{\partial}{\partial z}\otimes\frac{\partial}{\partial\overline{z}_{K}}\otimes\sigma^{*}\otimes e_{i_{1}}^{*}\wedge\cdots\wedge e_{i_{p}}^{*})\in T_{x}^{1,0}M$,
then we get from (47) that

\noindent $\left\Vert B_{\Phi}A\right\Vert ^{2}\,=\left\Vert B_{\Phi}A(\frac{\partial}{\partial z}\otimes\frac{\partial}{\partial\overline{z}_{K}}\otimes\sigma^{*}\otimes e_{i_{1}}^{*}\wedge\cdots\wedge e_{i_{p}}^{*})\right\Vert ^{2}$

\noindent $\ \ \ \ \ \ \ \ \ \ \ \ =\left\Vert B_{\Phi}(X_{\overline{K}I})\right\Vert ^{2}$

\noindent $\ \ \ \ \ \ \ \ \ \ \ \ =\left\Vert B_{\Phi}(X_{\overline{K}I})\cdot(e_{j_{1}}\wedge\cdots\wedge e_{j_{p}})\right\Vert ^{2}$

\noindent $\ \ \ \ \ \ \ \ \ \ \ \overset{(47)}{=}\left\Vert B_{s}(X_{\overline{K}I})\lrcorner(e_{j_{1}}\wedge\cdots\wedge e_{j_{p}})\right\Vert ^{2}$

\noindent $\ \ \ \ \ \ \ \ \ \ \ \ =\binom{r}{p-1}\left\Vert B_{s}(X_{\overline{K}I})\right\Vert ^{2}=\binom{r}{p-1}\left\Vert B_{s}A\right\Vert ^{2}.$

\noindent By (47) and the definition (12) , we have

\noindent $\textrm{Tr}B_{\Phi}A=(B_{\Phi}A(\frac{\partial}{\partial z}\otimes\frac{\partial}{\partial\overline{z}_{K}}\otimes\sigma^{*}\otimes e_{i_{1}}^{*}\wedge\cdots\wedge e_{i_{p}}^{*}),$

\noindent $\ \ \ \ \ \ \ \ \ \ \ \ \ \ e_{j_{1}}\wedge\cdots\wedge e_{j_{p-1}}\otimes e_{i_{1}}^{*}\wedge\cdots\wedge e_{i_{p}}^{*})dz\otimes\sigma\otimes e_{j_{1}}\wedge\cdots\wedge e_{j_{p-1}}$

\noindent $\ \ \ \ \ \ \ \ \ \ \ =(B_{\Phi}(X_{i_{1}\cdots i_{p}})\cdot(e_{i_{1}}\wedge\cdots\wedge e_{i_{p}}),e_{j_{1}}\wedge\cdots\wedge e_{j_{p-1}})$

\noindent $\ \ \ \ \ \ \ \ \ \ \ \ \ \ \cdot dz\otimes\sigma\otimes e_{j_{1}}\wedge\cdots\wedge e_{j_{p-1}}$

\noindent $\ \ \ \ \ \ \ \ \ \ \overset{(47)}{=}(B_{s}(X_{i_{1}\cdots i_{p}})\lrcorner(e_{i_{1}}\wedge\cdots\wedge e_{i_{p}}),e_{j_{1}}\wedge\cdots\wedge e_{j_{p-1}})$

\noindent $\ \ \ \ \ \ \ \ \ \ \ \ \ \ \cdot dz\otimes\sigma\otimes e_{j_{1}}\wedge\cdots\wedge e_{j_{p-1}}$

\noindent $\ \ \ \ \ \ \ \ \ \ \ =(B_{s}(X_{i_{1}\cdots i_{p}})\lrcorner(dz\otimes\sigma\otimes e_{i_{1}}\wedge\cdots\wedge e_{i_{p}})$

\noindent $\ \ \ \ \ \ \ \ \ \ \ =(B_{s}(A(\frac{\partial}{\partial z}\otimes\frac{\partial}{\partial\overline{z}_{K}}\otimes\sigma^{*}\otimes e_{i_{1}}^{*}\wedge\cdots\wedge e_{i_{p}}^{*}))$

\noindent $\ \ \ \ \ \ \ \ \ \ \ \ \ \ \lrcorner(dz\otimes\sigma\otimes e_{i_{1}}\wedge\cdots\wedge e_{i_{p}})$

\noindent $\ \ \ \ \ \ \ \ \ \ \ =\textrm{Tr}B_{s}A$.

\noindent We denote by $F,F^{p}$ the curvatures of $E$ and $\wedge^{p}E$
respectively($1\leq p\leq r)$. Let $\xi=\frac{1}{p!}\underset{1\leq i_{1},\cdots,i_{p}\leq r}{\sum}\xi_{i_{1}\cdots i_{p}}e_{i_{1}}\wedge\cdots\wedge e_{i_{p}},X=X_{\alpha}\frac{\partial}{\partial z_{\alpha}}.$
In all of the computations involving $\xi$ during the proof of the
equality (50), sum is always taken over all repeated indices(not only
over strictly increasing multi-indices). We know by definition that

\noindent $(F_{X\overline{X}}^{p}\xi,\xi)\,=\frac{1}{p!}(F_{X\overline{X}}^{p}\xi)_{i_{1}\cdots i_{p}}\overline{\xi_{i_{1}\cdots i_{p}}}$

\noindent $\ \ \ \ \ \ \ \ \ \ \ \ \ \ =\frac{1}{p!}\underset{1\leq a\leq p}{\sum}F_{\alpha\overline{\beta}i_{a}\overline{i}}X_{\alpha}\xi_{i_{1}\cdots(i)_{a}\cdots i_{p}}\overline{X_{\beta}\xi_{i_{1}\cdots i_{p}}}$

\noindent $\ \ \ \ \ \ \ \ \ \ \ \ \ \ =\frac{1}{p!}\underset{1\leq a\leq p}{\sum}F_{\alpha\overline{\beta}i_{a}\overline{i}}X_{\alpha}\xi_{ii_{1}\cdots\widehat{i_{a}}\cdots i_{p}}\overline{X_{\beta}\xi_{i_{a}i_{1}\cdots\widehat{i_{a}}\cdots i_{p}}}$

\noindent \begin{equation}
=\frac{1}{p!}\sum_{1\leq a\leq p}F_{\alpha\overline{\beta}i\overline{j}}X_{\alpha}\xi_{ji_{1}\cdots\widehat{i_{a}}\cdots i_{p}}\overline{X_{\beta}\xi_{ii_{1}\cdots\widehat{i_{a}}\cdots i_{p}}}.\ \ \ \ \ \ \ \ \ \ \ \ \ \ \ \ \ \end{equation}

\noindent where $\hat{i_{a}}$ means omitting the index $i_{a}.$
Similarly, we have

\noindent $(F_{X\overline{X}}^{p-1}\Phi\xi,\Phi\xi)=\frac{1}{(p-1)!}F_{\alpha\overline{\beta}i\overline{j}}X_{\alpha}(\Phi\xi)_{ji_{1}\cdots\widehat{i_{a}}\cdots i_{p-1}}\overline{X_{\beta}(\Phi\xi)_{ii_{1}\cdots\widehat{i_{a}}\cdots i_{p-1}}}$

\noindent $\,\ \ \ \ \ \ \ \ \ \ \ \ \ \ \ \ \ \ =\underset{1\leq a\leq p-1}{\sum}F_{\alpha\overline{\beta}i\overline{j}}X_{\alpha}g_{k}\xi_{jki_{1}\cdots\widehat{i_{a}}\cdots i_{p-1}}\overline{X_{\beta}g_{l}\xi_{ili_{1}\cdots\widehat{i_{a}}\cdots i_{p-1}}.}$

\noindent On the other hand, we also have

\noindent $(\Phi F_{X\overline{X}}^{p}\xi,\Phi\xi)=\frac{1}{(p-1)!}(\Phi F_{X\overline{X}}^{p}\xi)_{i_{1}\cdots i_{p-1}}\overline{(\Phi\xi)_{i_{1}\cdots i_{p-1}}}$

\noindent $\ \ \ \ \ \ \ \ \ \ \ \ \ \ \ \ \ \ =\frac{1}{(p-1)!}g_{k}(F_{X\overline{X}}^{p}\xi)_{ki_{1}\cdots i_{p-1}}\overline{g_{l}\xi_{li_{1}\cdots i_{p-1}}}$

\noindent $\ \ \ \ \ \ \ \ \ \ \ \ \ \ \ \ \ \ =\frac{1}{(p-1)!}F_{\alpha\overline{\beta}k\overline{i}}X_{\alpha}g_{k}\xi_{ii_{1}\cdots i_{p-1}}\overline{g_{l}\xi_{li_{1}\cdots i_{p-1}}}$

\noindent $\ \ \ \ \ \ \ \ \ \ \ \ \ \ \ \ \ \ \ \ \ +\frac{1}{(p-1)!}\underset{1\leq a\leq p-1}{\sum}F_{\alpha\overline{\beta}i_{a}\overline{j}}X_{\alpha}g_{k}\xi_{ki_{1}\cdots(j)_{a}\cdots i_{p-1}}\overline{X_{\beta}g_{l}\xi_{li_{1}\cdots i_{p-1}}}$

\noindent $\ \ \ \ \ \ \ \ \ \ \ \ \ \ \ \ \ \ =\frac{1}{(p-1)!}F_{\alpha\overline{\beta}k\overline{i}}X_{\alpha}g_{k}\xi_{ii_{1}\cdots i_{p-1}}\overline{X_{\beta}g_{l}\xi_{li_{1}\cdots i_{p-1}}}$

\noindent $\ \ \ \ \ \ \ \ \ \ \ \ \ \ \ \ \ \ \ \ \ +\frac{1}{(p-1)!}\underset{1\leq a\leq p-1}{\sum}F_{\alpha\overline{\beta}i_{a}\overline{j}}X_{\alpha}g_{k}\xi_{jki_{1}\cdots\widehat{i_{a}}\cdots i_{p-1}}\overline{X_{\beta}g_{l}\xi_{i_{a}li_{1}\cdots\widehat{i_{a}}\cdots i_{p-1}}}$

\noindent $\ \ \ \ \ \ \ \ \ \ \ \ \ \ \ \ \ \ =\frac{1}{(p-1)!}F_{\alpha\overline{\beta}k\overline{i}}X_{\alpha}g_{k}\xi_{ii_{1}\cdots i_{p-1}}\overline{X_{\beta}g_{l}\xi_{li_{1}\cdots i_{p-1}}}$

\noindent $\ \ \ \ \ \ \ \ \ \ \ \ \ \ \ \ \ \ \ \ \ +\frac{1}{(p-1)!}\underset{1\leq a\leq p-1}{\sum}F_{\alpha\overline{\beta}i\overline{j}}X_{\alpha}g_{k}\xi_{jki_{1}\cdots\widehat{i_{a}}\cdots i_{p-1}}\overline{X_{\beta}g_{l}\xi_{ili_{1}\cdots\widehat{i_{a}}\cdots i_{p-1}}}.$

\noindent Form the last two equalities, it follows that

\noindent $((F_{X\overline{X}}^{\textrm{Hom}(\wedge^{p}E,\wedge^{p-1}E)}\Phi)\xi,\Phi\xi)=(F_{X\overline{X}}^{p-1}\Phi\xi,\Phi\xi)-(\Phi F_{X\overline{X}}^{p}\xi,\Phi\xi)$

\noindent $\,\ \ \ \ \ \ \ \ \ \ \ \ \ \ \ \ \ \ \ \ \ \ \ \ \ \ \ \ \ \ \ \ \ \ \ \ =-\frac{1}{(p-1)!}F_{\alpha\overline{\beta}k\overline{i}}X_{\alpha}g_{k}\xi_{ii_{1}\cdots i_{p-1}}$

\noindent $\ \ \ \ \ \ \ \ \ \ \ \ \ \ \ \ \ \ \ \ \ \ \ \ \ \ \ \ \ \ \ \ \ \ \ \ \ \ \ \ \cdot\overline{X_{\beta}g_{l}\xi_{li_{1}\cdots i_{p-1}}}$

\noindent which implies

\noindent $(F_{X\overline{X}}^{\textrm{Hom}(\wedge^{p}E,\wedge^{p-1}E)}\Phi,\Phi)=\underset{i_{1}<\cdots<i_{p}}{\sum}(F_{X\overline{X}}^{\textrm{Hom}(\wedge^{p}E,\wedge^{p-1}E)}\Phi\cdot e_{i_{1}}\wedge\cdots\wedge e_{i_{p}}$

\noindent $\ \ \ \ \ \ \ \ \ \ \ \ \ \ \ \ \ \ \ \ \ \ \ \ \ \ \ \ \ \ \ \ \ \ \ ,\Phi\cdot e_{i_{1}}\wedge\cdots\wedge e_{i_{p}})$

\noindent $\ \ \ \ \ \ \ \ \ \ \ \ \ \ \ \ \ \ \ \ \ \ \ \ \ \ \ \ \ \ \ \ \ =-\underset{i_{1}<\cdots<i_{p-1}}{\sum}F_{\alpha\overline{\beta}k\overline{i}}X_{\alpha}g_{k}\overline{X_{\beta}g_{l}}sgn(\begin{array}{cccc}
i & i_{1} & \cdots & i_{p-1}\\
l & i_{1} & \cdots & i_{p-1}\end{array})$

\noindent $\ \ \ \ \ \ \ \ \ \ \ \ \ \ \ \ \ \ \ \ \ \ \ \ \ \ \ \ \ \ =-\binom{r}{p-1}F_{\alpha\overline{\beta}k\overline{l}}X_{\alpha}g_{k}\cdot\overline{X_{\beta}g_{l}}=\binom{r}{p-1}(F_{X\overline{X}}^{E^{*}}s,s).$

\begin{spacing}{0.7}
\noindent The proof is complete.
\end{spacing}

\begin{spacing}{0.7}
\noindent \begin{flushright}
$\square$$\ $
\par\end{flushright}
\end{spacing}

Consequently, we obtain the following identities.

$\ $

\noindent $\mathbf{Lemma5.2.}$ For any $u\in\wedge^{n,k-1}M\otimes L\otimes\wedge^{p-1}E$
and $v\in\wedge^{n,k}M\otimes L\otimes\wedge^{p}E,$ $1\leq k\leq n,$
we have the following pointwise identities outside $s^{-1}(0)$. \begin{equation}
\partial_{\alpha}\partial_{\overline{\beta}}\psi v_{\overline{\alpha K}I}\overline{v_{\overline{\beta K}I}}=e^{-\psi}\left\Vert B_{s}A_{v}\right\Vert ^{2}-e^{-\psi}(F_{X_{\overline{K}i}\overline{X_{\overline{K}i}}}^{E^{*}}s,s).\end{equation}

\begin{equation}
(\overline{\partial}\Phi^{*}\wedge u,v)-(\Phi^{*}u,\textrm{grad}^{0,1}\psi\lrcorner v)=(u,\textrm{Tr}B_{s}A_{v})\end{equation}

\noindent where $\psi=\textrm{log}|s|^{2},A_{v}$ is defined by (9),
$\Phi=s\lrcorner,$ $v=v_{\overline{K}I}dz\otimes d\overline{z}_{K}\otimes\sigma\otimes e_{I,}e_{I}=e_{i_{1}}\wedge\cdots\wedge e_{i_{p}},$
and $X_{\overline{K}i}=A_{v}(\frac{\partial}{\partial z}\otimes\frac{\partial}{\partial\overline{z}_{K}}\otimes\sigma^{*}\otimes e_{i}^{*}).$

$\,$

\noindent Proof. Let $\varphi=\log|\Phi|^{2},$ then we know, by (45),
$\varphi=\psi+\log$$\binom{r}{p-1}.$

\noindent From the identity (11), it follows that

\noindent $\partial_{\alpha}\partial_{\overline{\beta}}\psi v_{\overline{\alpha K}I}\overline{v_{\overline{\beta K}I}}=\partial_{\alpha}\partial_{\overline{\beta}}\varphi v_{\overline{\alpha K}I}\overline{v_{\overline{\beta K}I}}$

\noindent $\ \ \ \ \ \ \ \ \ \ \ \ \ \ \ \ \ \ \ \ \ \ =e^{-\varphi}[\left\Vert B_{\Phi}A_{v}\right\Vert ^{2}-(F_{X_{\overline{K}i}\overline{X_{\overline{K}i}}}^{\textrm{Hom}(E,E^{'})}\Phi,\Phi)]$

$\ \ \ \ \ \ \ \ \ \ \ \ \ \ \overset{(48),(50)}{=}e^{-\psi}\left\Vert B_{s}A_{v}\right\Vert ^{2}-e^{-\psi}(F_{X_{\overline{K}i}\overline{X_{\overline{K}i}}}^{E^{*}}s,s).$

\begin{spacing}{0.7}
\noindent (53) follows from (15) (45) (47). The proof is finished.
\end{spacing}

\begin{spacing}{0.7}
\noindent \begin{flushright}
$\ $$\ $$\square$
\par\end{flushright}
\end{spacing}

Now we improve, for $\Phi=s\lrcorner$, the main estimate obtained
in section 3.

$\ $

\noindent $\mathbf{Lemma5.3.}$ Let $(M,\omega)$ be a K\"{a}hler
manifold and let $E$ be a Hermitian holomorphic vector bundle over
$M$, $L$ a Hermitian holomorphic line bundle over $M$.
$\Omega\Subset M\setminus s^{-1}(0)$ is a pseudoconvex domain with
smooth boundary. Suppose that the following conditions hold on
$\Omega.$

1. $E\geq_{m}0,m\geq\textrm{min}\{n-k+1,r-p+1\},1\leq k\leq n,1\leq p\leq r;$

2. the curvature of $L$ satisfies

\begin{center}
$\sqrt{-1}(\varsigma c(L)-\partial\overline{\partial}\varsigma-\tau^{-1}\partial\varsigma\wedge\overline{\partial}\varsigma)\geq\sqrt{-1}q(\varsigma+\delta)\partial\overline{\partial}\psi$.
\par\end{center}

\noindent Then the following estimate \begin{equation}
\left\Vert |s|^{-2}\theta\wedge u+\overline{\partial}^{*}v\right\Vert _{\Omega,\varsigma+\tau}^{2}+\left\Vert \overline{\partial}v\right\Vert _{\Omega,\varsigma}^{2}\geq\left\Vert u\right\Vert _{\Omega,\frac{\varsigma\delta}{(\varsigma+\delta)|s|^{2}}}^{2}\end{equation}
holds for every $\overline{\partial}$-closed $u\in A^{n,k-1}(\overline{\Omega},L\otimes\wedge^{p-1}E)$
satisfying $s\lrcorner u=0$ and every $v\in A^{n,k}(\overline{\Omega},L\otimes\wedge^{p}E)\cap\textrm{Dom}(\overline{\partial}^{*}),$
where $\theta$ is defined in (44), $\psi=\log\left|s\right|^{2},q=\textrm{min}\{n,r-1\},n=\dim_{\mathbb{C}}M,r=\textrm{rank}_{\mathbb{C}}E,$
$0<\varsigma\in C^{\infty}(\Omega)$ and $\delta,\tau$ are measurable
functions on $\Omega$ satisfying $\tau>0,\varsigma+\delta\geq0.$

$\ $

\noindent Proof. The proof is essentially the same as that of lemma
3.3, we sketch it with an emphasis on the modifications.

\noindent By the following identity

\noindent \begin{center}
$|\Phi^{*}u|^{2}=(\theta\wedge u,\theta\wedge u)=(s\lrcorner\theta\wedge u,u)$
\par\end{center}

\noindent \begin{center}
$\ \ \ \ \ \ \ $$\ \ \ \ \ \ \ \ \ \ \ \ \ \ \ \ \ \ \ \ \ \ \ \ =(s\lrcorner\theta\wedge u+\theta\wedge s\lrcorner u,u)=(|s|^{2}u,u)=|s|^{2}|u|^{2},$
\par\end{center}

\noindent we obtain

\noindent l.h.s.of (54)$=\left\Vert |s|^{-2}\Phi^{*}u+\overline{\partial}^{*}v\right\Vert _{\varsigma+\tau}^{2}+\left\Vert \overline{\partial}v\right\Vert _{\varsigma}^{2}$

\noindent $\ \ \ \ \ \ \ \ \ \ \ \ \ \ \ =\left\Vert u\right\Vert _{|s|^{-2}\varsigma}^{2}+\left\Vert |s|^{-2}\Phi^{*}u+\overline{\partial}^{*}v\right\Vert _{\tau}^{2}$

\noindent $\ \ \ \ \ \ \ \ \ \ \ \ \ \ \ \ \ \ +2\textrm{Re}(\varsigma e^{-\psi}\Phi^{*}u,\overline{\partial}^{*}v)+\left\Vert \sqrt{\varsigma}\overline{\partial}^{*}v\right\Vert ^{2}+\left\Vert \sqrt{\varsigma}\overline{\partial}v\right\Vert ^{2}.$

\noindent By direct computation, we have

\noindent $([\sqrt{-1}F^{L\otimes\wedge^{p}E},\Lambda_{\omega}]v,v)=F_{\alpha\overline{\beta}I\overline{J}}^{L\otimes\wedge^{p}E}v_{\overline{\alpha K},I}\overline{v_{\overline{\beta K},J}}$

$\ \ \ \ \ \ \ \ \ \ \ \ \ \ \ \ \ \ \ \ \ \ \ \ \ \ \ =F_{\alpha\overline{\beta}I\overline{J}}^{\wedge^{p}E}v_{\overline{\alpha K},I}\overline{v_{\overline{\beta K},J}}+F_{\alpha\overline{\beta}}^{L}v_{\overline{\alpha K},I}\overline{v_{\overline{\beta K},I}}$

$\ \ \ \ \ \ \ \ \ \ \ \ \ \ \ \ \ \ \ \ \ \ \ \ \ \ \ =F_{\alpha\overline{\beta}i\overline{j}}^{E}v_{\overline{\alpha K},iN}\overline{v_{\overline{\beta K},iN}}+F_{\alpha\overline{\beta}}^{L}v_{\overline{\alpha K},I}\overline{v_{\overline{\beta K},I}}$

$\ \ \ \ \ \ \ \ \ \ \ \ \ \ \ \ \ \ \ \ \ \ \ \ \ \ \ \geq F_{\alpha\overline{\beta}}^{L}v_{\overline{\alpha K},I}\overline{v_{\overline{\beta K},I}},$

\noindent the last inequality follows from the condition $E\geq_{m}0,m\geq\textrm{min}\{n-k+1,r-p+1\}.$
Now we get by using the twisted Bochner-Kodaira-Nakano formula and
Morrey's trick(to handle the boundary term)

\noindent l.h.s.of (54)$\geq\left\Vert u\right\Vert _{|s|^{-2}\varsigma}^{2}+\left\Vert |s|^{-2}\Phi^{*}u+\overline{\partial}^{*}v\right\Vert _{\tau}^{2}$

\noindent $\ \ \ \ \ \ \ \ \ \ \ \ \ \ \ \ \ \ +\int_{\Omega}\varsigma F_{\alpha\overline{\beta}}^{L}v_{\overline{\alpha K},I}\overline{v_{\overline{\beta K},I}}-\nabla^{\overline{\alpha}}\nabla^{\beta}\varsigma(\frac{\partial}{\partial\overline{z}_{\alpha}}\lrcorner v,\frac{\partial}{\partial\overline{z}_{\beta}}\lrcorner v)$

\noindent $\ \ \ \ \ \ \ \ \ \ \ \ \ \ \ \ \ \ +2\textrm{Re}(\varsigma e^{-\psi}\Phi^{*}u+\textrm{grad}^{0,1}\varsigma\lrcorner v,\overline{\partial}^{*}v)dV_{\omega}$

\noindent $\ \ \ \ \ \ \ \ \ \ \ \ \ \ \ \overset{(51)}{\geq}\left\Vert u\right\Vert _{|s|^{-2}\varsigma}^{2}+\left\Vert |s|^{-2}\Phi^{*}u+\overline{\partial}^{*}v\right\Vert _{\tau}^{2}$

\noindent $\ \ \ \ \ \ \ \ \ \ \ \ \ \ \ \ \ \ +\int_{\Omega}\varsigma F_{\alpha\overline{\beta}}^{L}v_{\overline{\alpha K},I}\overline{v_{\overline{\beta K},I}}-\nabla^{\overline{\alpha}}\nabla^{\beta}\varsigma(\frac{\partial}{\partial\overline{z}_{\alpha}}\lrcorner v,\frac{\partial}{\partial\overline{z}_{\beta}}\lrcorner v)$

\noindent $\ \ \ \ \ \ \ \ \ \ \ \ \ \ \ \ \ \ +2\textrm{Re}(\varsigma e^{-\psi}\Phi^{*}u+\textrm{grad}^{0,1}\varsigma\lrcorner v,\overline{\partial}^{*}v)dV_{\omega}$

\noindent $\ \ \ \ \ \ \ \ \ \ \ \ \ \ \ \geq\left\Vert u\right\Vert _{|s|^{-2}\varsigma}^{2}+\left\Vert |s|^{-2}\Phi^{*}u+\overline{\partial}^{*}v\right\Vert _{\tau}^{2}$

\noindent $\ \ \ \ \ \ \ \ \ \ \ \ \ \ \ \ \ \ +\int_{\Omega}(q(\varsigma+\delta)\partial_{\alpha}\partial_{\overline{\beta}}\psi+\tau^{-1}\partial_{\alpha}\varsigma\overline{\partial_{\beta}\varsigma})v_{\overline{\alpha K},I}\overline{v_{\overline{\beta K},I}}$

\noindent $\ \ \ \ \ \ \ \ \ \ \ \ \ \ \ \ \ \ +2\textrm{Re}(\varsigma e^{-\psi}\Phi^{*}u+\textrm{grad}^{0,1}\varsigma\lrcorner v,\overline{\partial}^{*}v)dV_{\omega}$

\noindent $\ \ \ \ \ \ \ \ \ \ \ \ \ \ \ \overset{(52)}{\geq}\left\Vert u\right\Vert _{|s|^{-2}\varsigma}^{2}+\left\Vert |s|^{-2}\Phi^{*}u+\overline{\partial}^{*}v\right\Vert _{\tau}^{2}+\left\Vert \sqrt{q(\varsigma+\delta)}B_{s}A_{v}\right\Vert _{|s|^{-2}}^{2}$

\noindent \begin{equation}
\ \ \ \ \ \ \ +\left\Vert \textrm{grad}^{0,1}\varsigma\lrcorner v\right\Vert _{\tau^{-1}}^{2}+2\textrm{Re}(\varsigma e^{-\psi}\Phi^{*}u+\textrm{grad}^{0,1}\varsigma\lrcorner v,\overline{\partial}^{*}v)\end{equation}

\noindent Integration by parts yields

\noindent $2\textrm{Re}(\varsigma e^{-\psi}\Phi^{*}u,\overline{\partial}^{*}v)=2\textrm{Re}(e^{-\psi}\Phi^{*}u,\textrm{grad}^{0,1}\varsigma\lrcorner v)+2\textrm{Re}(\overline{\partial}\Phi^{*}\wedge u,v)_{|s|^{-2}\varsigma}$

\noindent $\ \ \ \ \ \ \ \ \ \ \ \ \ \ \ \ \ \ \ \ \ \ \ \ \ \ \ \ \ -2\textrm{Re}(\Phi^{*}u,\textrm{grad}^{0,1}\psi\lrcorner v)_{|s|^{-2}\varsigma}$

\noindent $\ \ \ \ \ \ \ \ \ \ \ \ \ \ \ \ \ \ \ \ \ \ \ \ \ \ \overset{(53)}{=}2\textrm{Re}(e^{-\psi}\Phi^{*}u,\textrm{grad}^{0,1}\varsigma\lrcorner v)+2\textrm{Re}(u,\textrm{Tr}B_{s}A_{v})_{\left|s\right|^{-2}\varsigma}$

\noindent $\ \ \ \ \ \ \ \ \ \ \ \ \ \ \ \ \ \ \ \ \ \ \ \ \ \ \ \geq2\textrm{Re}(e^{-\psi}\Phi^{*}u,\textrm{grad}^{0,1}\varsigma\lrcorner v)$

\noindent \begin{equation}
\ \ \ \ \ \ \ \ \ \ \ -\left\Vert \sqrt{\varsigma+\delta}\textrm{Tr}B_{s}A_{v}\right\Vert _{|s|^{-2}}^{2}-\left\Vert u\right\Vert _{\frac{\varsigma^{2}}{(\varsigma+\delta)|s|^{2}}}^{2}.\end{equation}

\noindent By definition, we have

\begin{center}
$\textrm{rank}B_{s}A_{v}\leq\textrm{rank}B_{s}\leq\textrm{min}\{n,r-1\}.$
\par\end{center}

\noindent It follows from (13) (55) and (56) that

\noindent l.h.s.of (55)$\geq\left\Vert u\right\Vert _{\frac{\varsigma\delta}{(\varsigma+\delta)|s|^{2}}}^{2}+\left\Vert |s|^{-2}\Phi^{*}u+\overline{\partial}^{*}v\right\Vert _{\tau}^{2}$

\noindent $\ \ \ \ \ \ \ \ \ \ \ \ \ \ \ \ \ \ +\left\Vert \textrm{grad}^{0,1}\varsigma\lrcorner v\right\Vert _{\tau^{-1}}^{2}+2\textrm{Re}(\textrm{grad}^{0,1}\varsigma\lrcorner v,e^{-\psi}\Phi^{*}u+\overline{\partial}^{*}v)$

\noindent $\ \ \ \ \ \ \ \ \ \ \ \ \ \ \ \geq\left\Vert u\right\Vert _{\frac{\varsigma\delta}{(\varsigma+\delta)|s|^{2}},}^{2}$

\begin{spacing}{0.7}
\noindent which completes the proof.
\end{spacing}

\begin{spacing}{0.7}
\noindent \begin{flushright}
$\square$
\par\end{flushright}
\end{spacing}

$\ $By the standard functional argument used in the proof of theorem
4.2 with the estimate (20) replaced by the improved one (54), we obtain
the main result of this section:

$\ $

\noindent $\mathbf{Theorem5.4.}$ Let $(M,\omega)$ be a
K$\ddot{a}$ler manifold and let $E$ be a Hermitian holomorphic
vector bundle over $M$, $L$ a line bundle over $M,$ $s\in$
$\Gamma(E^{*})$. All the Hermitian structures may have singularities
in a subvariety $Z\subsetneqq M$ $.$ We define the Koszul complex
associated to $s$ by (42) (43). Assume that $s^{-1}(0)\subseteq Z,$
and that $M\setminus Z$ is weakly pseudoconvex and that the
following conditions hold on $M\setminus Z$:

1. $E\geq_{m}0,m\geq\textrm{min}\{n-k+1,r-p+1\},1\leq k\leq n,1\leq p\leq n;$

2. the curvature of $L$ satisfies

\begin{center}
$\sqrt{-1}(\varsigma c(L)-\partial\overline{\partial}\varsigma-\tau^{-1}\partial\varsigma\wedge\overline{\partial}\varsigma)\geq\sqrt{-1}q(\varsigma+\delta)\partial\overline{\partial}\varphi$.
\par\end{center}

\noindent Then for any $\overline{\partial}$-closed $(n,k-1)$-form
$f$ which is valued in $L\otimes\wedge^{p-1}E,$ if $d_{p-1}f=0$
and $\left\Vert f\right\Vert _{\frac{\varsigma+\delta}{\varsigma\delta|s|^{2}}}<+\infty$
then there is at least one $\overline{\partial}$-closed $(n,k-1)$-form
$h$ valued in $L\otimes\wedge^{p}E$ such that $d_{p}h=f$ and the
following estimate holds \begin{equation}
\left\Vert h\right\Vert _{\frac{1}{\varsigma+\tau}}\leq\left\Vert f\right\Vert _{\frac{\varsigma+\delta}{\varsigma\delta|s|^{2}}},\end{equation}
 where $1\leq p\leq r,\varphi=\log\left|s\right|,q=\textrm{min}\{n,r-1\},n=\dim_{\mathbb{C}}M,r=\textrm{rank}_{\mathbb{C}}E,$
$0<\varsigma,\tau\in C^{\infty}(M)$ and $\delta\geq0$ is a
measurable function on $M.$\\

We can derive from theorem 5.4 the next result by repeating the
argument used in the proof of corollary 4.3.\\

\noindent $\mathbf{Corollary5.5.}$ Let $(M,\omega)$ be a K\"{a}hler
manifold and let $E$ be a Hermitian holomorphic vector bundle over
$M$, $L$ a line bundle over $M,$ $s\in$ $\Gamma(E^{*})$. All the
Hermitian structures may have singularities in a subvariety
$Z\subsetneqq M$ $.$ We define the Koszul complex associated to $s$
by (42) (43). Assume that $s^{-1}(0)\subseteq Z,$ and that
$M\setminus Z$ is weakly pseudoconvex and the following conditions
hold on $M\setminus Z$:

1. $E\geq_{m}0,m\geq\textrm{min}\{n-k+1,r-p+1\},1\leq k\leq n,1\leq p\leq r;$

2. the curvature of $L$ satisfies $\sqrt{-1}c(L)\geq\sqrt{-1}q(1+\varepsilon)\partial\overline{\partial}\varphi$.

\noindent Then for any $\overline{\partial}$-closed $(n,k-1)$-form
$f$ valued in $L\otimes\wedge^{p-1}E,$ if $d_{p-1}f=0$ and $\left\Vert f\right\Vert _{|s|^{-2}}<+\infty$
then there is at least one $\overline{\partial}$-closed $(n,k-1)$-form
$h$ valued in $L\otimes\wedge^{p}E$ such that $d_{p}h=f$ with the
estimate \begin{equation}
\left\Vert h\right\Vert ^{2}\leq\frac{1+\varepsilon}{\varepsilon}\left\Vert f\right\Vert _{|s|^{-2}}^{2},\end{equation}
 where $1\leq p\leq r,\varphi=\log\left|s\right|^{2},q=\textrm{min}\{n,r-1\},n=\dim_{\mathbb{C}}M,r=\textrm{rank}_{\mathbb{C}}E$
and $\varepsilon$ is a positive constant.\\

\noindent As a consequence of corollary 4.5, we also have the
following sufficient condition for the exactness of the induced
sequence of global sections.\\

\noindent $\mathbf{Corollary5.6.}$ Let $(M,\omega)$ be a weakly
pseudoconvex K\"{a}hler manifold and let $E$ be a Hermitian
holomorphic vector bundle over $M$, $L$ a line bundle over $M$.
Assume that $s\in$ $\Gamma(E^{*})$ is a nonvanishing section and
that

1. $E$ is semi-positive in the sense of Nakano;

2. the curvature of $L$ satisfies

\begin{center}
$\sqrt{-1}c(L)\geq\sqrt{-1}q(1+\varepsilon)\partial\overline{\partial}\varphi$
for some positive constant $\varepsilon$.
\par\end{center}

\noindent Then the induced sequence on global sections

\begin{center}
$\Gamma(K_{M}\otimes L\otimes\textrm{det}E)\overset{d_{r}}{\rightarrow}\Gamma(K_{M}\otimes L\otimes\wedge^{r-1}E)\overset{d_{r-1}}{\rightarrow}\cdots\overset{d_{1}}{\rightarrow}\Gamma(K_{M}\otimes L)\overset{d_{0}}{\rightarrow}0$
\par\end{center}

\noindent is exact.

$\,$

Now we discuss the special case of Koszul complex over a domain
$\Omega\subseteq\mathbb{C}^{n}.$

Let $g_{1}\cdots,g_{r}\in\mathcal{O}(\Omega)$, the Koszul complex
associated to $g=(g_{1}\cdots,g_{r})$ is given by \begin{equation}
\wedge^{r}\mathcal{O}^{\oplus r}\overset{d_{r}}{\rightarrow}\wedge^{r-1}\mathcal{O}^{\oplus r}\overset{d_{r-1}}{\rightarrow}\cdots\overset{d_{2}}{\rightarrow}\wedge\mathcal{O}^{\oplus r}\overset{d_{1}}{\rightarrow}\mathcal{O}\overset{d_{0}}{\rightarrow}0\end{equation}

\noindent where the boundary operators are defined by $d_{p}=g\lrcorner,1\leq p\leq r.$

\noindent It is easy to see that for every $h=(h_{i_{1}\cdots i_{p}})_{i_{1}\cdots i_{p}=1}^{r}\in\Gamma(\Omega,\wedge^{p}\mathcal{O}^{\oplus r})$$(\textrm{i.e. }h_{i_{1}\cdots i_{p}}\in\mathcal{O}(\Omega)$
and $h_{i_{1}\cdots i_{p}}$ is skew symmetric in $i_{1},\cdots,i_{p}),$we
have

\begin{center}
$d_{p}h=(f_{i_{1}\cdots i_{p-1}})_{i_{1}\cdots i_{p-1}=1}^{r}\textrm{\ensuremath{\in\Gamma}(\ensuremath{\Omega},\ensuremath{\wedge^{p-1}\mathcal{O}^{\oplus r}}) with}\ f_{i_{1}\cdots i_{p-1}}=\underset{1\leq\nu\leq r}{\sum}g_{\nu}h_{\nu i_{1}\cdots i_{p-1}}.$
\par\end{center}

Given a measurable function $\phi$ on $\Omega$ which is locally
bounded from above, we can define the following space

\begin{center}
$\left\{ (h_{i_{1}\cdots i_{p}})_{i_{1}\cdots i_{p}=1}^{r}\in\Gamma(\Omega,\wedge^{p}\mathcal{O}^{\oplus r})\mid\underset{i_{1}<\cdots<i_{p}}{\sum}\int_{\Omega}|h|^{2}e^{-\phi}dV<+\infty\right\} $
\par\end{center}

\noindent where
$\left|h\right|^{2}:=\underset{i_{1}<\cdots<i_{p}}{\sum}|h_{i_{1}\cdots
i_{p}}|^{2}.$ This space will be denoted by
$A_{\Omega}^{2}(\wedge^{p}\mathcal{O}^{\oplus r},\phi)$, $1\leq
p\leq r$. Since the measurable function $\phi$ is locally bounded
from above, we know that
$A_{\Omega}^{2}(\wedge^{p}\mathcal{O}^{\oplus r},\phi)$ is a Hilbert
space.\\

\noindent As a consequence of corollary 5.5, it follows that\\

\noindent $\mathbf{Corollary5.7.}$ Suppose $\Omega$ is a
pseudoconvex domain
in\emph{$\mathbb{C}^{n},\psi\in\rm{PSH}(\Omega),$
}$g_{1}\cdots,g_{r}\in\mathcal{O}(\Omega)$, $1\leq p\leq r$, and
$\varepsilon>0$ is a constant. For every cycle $f$ of the Koszul
complex of holomorphic functions, if $f\in
A_{\Omega}^{2}(\wedge^{p-1}\mathcal{O}^{\oplus
r},\psi+(q+q\varepsilon+1)\log\left|g\right|)$ , then there exists
at least one holomorphic $h\in$
$A_{\Omega}^{2}(\wedge^{p}\mathcal{O}^{\oplus
r},\psi+q(1+\varepsilon)\textrm{log}\left|g\right|^{2})$ such that
$f$ is the image of $h$ under the boundary map and

\begin{equation}
\int_{\Omega}|h|^{2}|g|^{-2q(1+\varepsilon)}e^{-\psi}dV\leq\frac{1+\varepsilon}{\varepsilon}\int_{\Omega}|f|^{2}|g|^{-2(q+q\varepsilon+1)}e^{-\psi}dV.\end{equation}

\noindent where $\left|g\right|=\sqrt{\sum_{i}\left|g_{i}\right|^{2}},q=\textrm{min}\{n,r-1\}.$

Particularly, if $\left|g\right|\neq0$ holds on $\Omega$ then (59)
induces an exact sequence on global sections:

\begin{center}
$\Gamma(\Omega,\wedge^{r}\mathcal{O}^{\oplus r})\overset{d_{r}}{\rightarrow}\Gamma(\Omega,\wedge^{r-1}\mathcal{O}^{\oplus r})\overset{d_{r-1}}{\rightarrow}\cdots\overset{d_{2}}{\rightarrow}\Gamma(\Omega,\wedge\mathcal{O}^{\oplus r})\overset{d_{1}}{\rightarrow}\Gamma(\Omega,\mathcal{O})\overset{d_{0}}{\rightarrow}0.$
\par\end{center}

$\:$

\noindent Proof. Step 1. We first give a proof in the case of
\emph{$\psi\in\rm{PSH}(\Omega)\cap C^{\infty}(\Omega).$ }

\noindent Let $E=\mathcal{O}^{\oplus r},L=\mathcal{O},s=(g_{1},\cdots,g_{r})\in\Gamma(E^{*}).$
We define the Hermitian structure on $E$ to be $h_{0},$ the standard
Hermitian structure on $\mathcal{O}^{\oplus r}$ induced form $\mathbb{C}^{r},$
so $(E,h_{0})$ is flat.

\noindent Let $\varepsilon>0$, then the following function\[
\frac{1}{(\sum_{i}\left|g_{i}\right|^{2})^{q(1+\varepsilon)}e^{\psi}}\]

\noindent defines a Hermitian structure on $L$ which has singularity
in $Z:=g_{1}^{-1}(0)$ (without loss of generality, we assume $g_{1}$
is not identically zero). The curvature this Hermitian structure is
given by

$\sqrt{-1}c(L)=\sqrt{-1}\overline{\partial}\partial(\textrm{log}(\sum_{i}\left|g_{i}\right|^{2})^{-q(1+\varepsilon)}e^{-\psi})$

$\ \ \ \ \ \ \ \ \ \ \ \ \ =q(1+\varepsilon)\sqrt{-1}\partial\overline{\partial}\textrm{log}|s|^{2}+\sqrt{-1}\partial\overline{\partial}\psi$

$\ \ \ \ \ \ \ \ \ \ \ \ \ \geq q(1+\varepsilon)\sqrt{-1}\partial\overline{\partial}\textrm{log}|s|^{2}.$

\noindent So far we have checked that all of the conditions assumed
in corollary 5.5 are fulfilled by the Hermitian structures on $E,L$
and $s\in\Gamma(E^{*})$ as constructed above. Hence the desired solvability
and estimate in this case follow from corollary 5.5.

Step 2. Now we proceed to prove the general case where $\psi$ is
only assumed to be plurisubharmonic on $\Omega.$

\noindent Let $\Omega_{1}\Subset\Omega_{2}\Subset\cdots$ be a
pseudoconvex exhaustion of $\Omega,$ and let
$\psi_{\ell}\in\rm{PSH}(\Omega_{\ell})\cap
C^{\infty}(\Omega_{\ell})\;\textrm{for }\ell\geq1$ such that

\begin{center}
$\psi_{\ell}\searrow\psi$ as $\ell\nearrow\infty$ on every $\Omega_{j}.$
\par\end{center}

\noindent Fix some $j\geq1,$ without loss of generality, we could
assume $\psi_{\ell}\in\rm{PSH}(\Omega_{j})\cap
C^{\infty}(\Omega_{j}),\ell\geq1$. Since $\psi_{\ell}\geq\psi$ holds
on $\Omega_{j},$ we get

\begin{center}
$\int_{\Omega_{j}}|f|^{2}|g|^{-2(q+q\varepsilon+1)}e^{-\psi_{\ell}}dV\leq\int_{\Omega}|f|^{2}|g|^{-2(q+q\varepsilon+1)}e^{-\psi}dV<+\infty$
\par\end{center}

\noindent which implies(by the result proved in step 1) that there
exists a holomorphic $h_{j\ell}\in$ $A_{\Omega_{j}}^{2}(\wedge^{p}\mathcal{O}^{\oplus r},\psi_{\ell}+q(1+\varepsilon)\textrm{log}\left|g\right|^{2})$
such that

\noindent \begin{center}
$d_{p}h_{j\ell}=f\mid_{\Omega_{j}}$
\par\end{center}

\noindent and

\begin{center}
$\int_{\Omega_{j}}|h_{j\ell}|^{2}|g|^{-2q(1+\varepsilon)}e^{-\psi_{1}}dV\leq\int_{\Omega_{j}}|h_{j\ell}|^{2}|g|^{-2q(1+\varepsilon)}e^{-\psi_{\ell}}dV$
\par\end{center}

\begin{center}
$\ \ \ \ \ \ \ \ \ \ \ \ \ \ \ \ \ \ \ \ \ \ \ \ \ \ \ \ \ \ \ \ \ \ \ \ \ \ \ \ \ \ \ \leq\frac{1+\varepsilon}{\varepsilon}\int_{\Omega_{j}}|f|^{2}|g|^{-2(q+q\varepsilon+1)}e^{-\psi_{\ell}}dV$
\par\end{center}

\begin{center}
$\ \ \ \ \ \ \ \ \ \ \ \ \ \ \ \ \ \ \ \ \ \ \ \ \ \ \ \ \ \ \ \ \ \ \ \ \ \ \ \ \ \ \leq\frac{1+\varepsilon}{\varepsilon}\int_{\Omega}|f|^{2}|g|^{-2(q+q\varepsilon+1)}e^{-\psi}dV.$
\par\end{center}

\noindent By the above estimate, we know $\{h_{j\ell}\}$ is bounded
in $A_{\Omega_{j}}^{2}(\wedge^{p}\mathcal{O}^{\oplus r},\psi_{1}+q(1+\varepsilon)\log\left|g\right|)$.
Consequently we can find a weak limit $h_{j}$ of $\{h_{j\ell}\}$
as $\ell\rightarrow\infty$ in $A_{\Omega_{j}}^{2}(\wedge^{p}\mathcal{O}^{\oplus r},\psi_{1}+q(1+\varepsilon)\log\left|g\right|).$

\noindent It is easy to see that

\noindent \begin{center}
$d_{p}h_{j}=f\mid_{\Omega_{j}}$
\par\end{center}

\noindent and

\noindent \begin{center}
$\ \ $$\int_{\Omega_{j}}|h_{j}|^{2}|g|^{-2q(1+\varepsilon)}e^{-\psi_{1}}dV\leq\underset{\ell\rightarrow\infty}{\underline{\textrm{lim}}}\int_{\Omega_{j}}|h_{j\ell}|^{2}|g|^{-2q(1+\varepsilon)}e^{-\psi_{1}}dV$
\par\end{center}

\begin{center}
$\ \ \ \ \ \ \ \ \ \ \ \ \ \ \ \ \ \ \ \ \ \ \ \ \ \ \ \ \ \ \ \ \ \ \ \ \leq\frac{1+\varepsilon}{\varepsilon}\int_{\Omega}|f|^{2}|g|^{-2(q+q\varepsilon+1)}e^{-\psi}dV.$
\par\end{center}

\begin{spacing}{0.7}
\noindent From the resulting estimate on $h_{j},$ we can repeat the
same argument of taking a weak limit of $\{h_{j}\}$ when $j\rightarrow\infty$
and then use the standard Cantor diagonalization process to show the
existence of the desired section $h$, it finishes our proof of corollary
5.7.
\end{spacing}

\begin{spacing}{0.7}
\noindent \begin{flushright}
$\square$
\par\end{flushright}
\end{spacing}

$\:$The special case of $p=1$ in corollary 5.7 is Skoda's division
theorem.\\

\noindent $\mathbf{Corollary5.8.}$ Suppose $\Omega$ is a
pseudoconvex domain in\emph{
$\mathbb{C}^{n},\psi\in\rm{PSH}(\Omega),$
}$g_{1}\cdots,g_{r}\in\mathcal{O}(\Omega)$, and $\varepsilon>0$ is a
constant, then for every $f\in
A_{\Omega}^{2}(\psi+(q+q\varepsilon+1)\log\left|g\right|^{2})$,
there exist holomorphic functions $h_{1},\cdots,h_{r}\in
A_{\Omega}^{2}(\psi+q(1+\varepsilon)\log\left|g\right|^{2})$ such
that $f=\sum_{i}g_{i}h_{i}$ and

\begin{equation}
\int_{\Omega}|h|^{2}|g|^{-2q(1+\varepsilon)}e^{-\psi}dV\leq\frac{1+\varepsilon}{\varepsilon}\int_{\Omega}|f|^{2}|g|^{-2(q+q\varepsilon+1)}e^{-\psi}dV.\end{equation}

\noindent
where$\left|g\right|^{2}=\sum_{i}\left|g_{i}\right|^{2},\left|h\right|^{2}=\sum_{i}\left|h_{i}\right|^{2},q=\textrm{min}\{n,r-1\}$
.\\

\noindent We end up this section by giving a sufficient condition
which is deduced from Skoda's division theorem by purely algebraic
argument.\\

\noindent Let $\Omega$ be a domain in $\mathbb{C}^{n},$ and $\Phi$
be a $q\times p$ matrix of holomorphic functions on $\Omega,p\geq q.$
We denote by $\delta_{i_{1}\cdots i_{q}}$ the $q\times q$ minors
of $\Phi,$ i.e.

\begin{center}
$\delta_{i_{1}\cdots i_{q}}=det\left(\begin{array}{ccc}
\Phi_{1i_{1}} & \cdots & \Phi_{1i_{q}}\\
\vdots & \ddots & \vdots\\
\Phi_{qi_{1}} & \cdots & \Phi_{qi_{q}}\end{array}\right),$
\par\end{center}

\noindent where $1\leq i_{1}<i_{2}<\cdots<i_{q}\leq p.$ There are
$\binom{p}{q}$ distinct minors of order $q.$

$\ $

\noindent $\mathit{\mathbf{Propsition5.9.}}$ Let
$\psi\in\rm{PSH}(\Omega),$ $f\in\mathcal{O}^{q}(\Omega),$ if
$\Omega$ is pseudoconvex and there exists a constant $\alpha>1$ such
that

\begin{center}
\begin{equation}
\int_{\Omega}\frac{|f|^{2}}{(\underset{i_{1}<\cdots<i_{q}}{\sum}|\delta_{i_{1}\cdots i_{q}}|^{2})^{\beta}}e^{-\psi}dV<+\infty,\end{equation}

\par\end{center}

\noindent where
$\beta=\textrm{min}\{n,\binom{p}{q}-1\}\cdot\alpha+1.$ Then there is
at least one $h\in\mathcal{O}^{p}(\Omega)$ which solves the
equations $\Phi h=f.$\\

 \noindent Proof. Let $f_{1},\cdots,f_{q}$ be the components of
$f,$ since for each $1\leq\nu\leq q$ we have

\begin{center}
$\ $$\int_{\Omega}\frac{|f_{\nu}|^{2}}{(\underset{i_{1}<\cdots<i_{q}}{\sum}|\delta_{i_{1}\cdots i_{q}}|^{2})^{\beta}}e^{-\psi}dV\leq\int_{\Omega}\frac{|f|^{2}}{(\underset{i_{1}<\cdots<i_{q}}{\sum}|\delta_{i_{1}\cdots i_{q}}|^{2})^{\beta}}e^{-\psi}dV<+\infty,$
\par\end{center}

\noindent there exists, by Skoda's theorem(i.e. corollary 5.8), a
system of functions $u_{i_{1}\cdots i_{q},\nu}\mathcal{\in O}^{q}(\Omega),1\leq i_{1}<i_{2}<\cdots<i_{q}\leq p,$
such that

\noindent \begin{center}
$f_{\nu}=\underset{i_{1}<\cdots<i_{q}}{\sum}\delta_{i_{1}\cdots i_{q}}u_{i_{1}\cdots i_{q},\nu}.$

\par\end{center}

\noindent Set

\noindent \begin{center}
$u_{i_{1}\cdots i_{q}}=\left(\begin{array}{c}
u_{i_{1}\cdots i_{q},1}\\
\vdots\\
u_{i_{1}\cdots i_{q},q}\end{array}\right),$
\par\end{center}

\noindent then the above equality could be rewritten as \begin{equation}
f=\underset{i_{1}<\cdots<i_{q}}{\sum}\delta_{i_{1}\cdots i_{q}}u_{i_{1}\cdots i_{q}}.\end{equation}

\begin{singlespace}
\noindent For fixed $1\leq i_{1}<i_{2}<\cdots<i_{q}\leq p,$ we consider
the $p$ by $p$ matrix $\Phi_{i_{1}\cdots i_{q}}$ whose entries
are defined by

\noindent \begin{equation}
\Phi_{i_{1}\cdots i_{q},ij}=\begin{cases}
\begin{array}{ccc}
\Phi_{ij}, & \textrm{if} & 1\leq i\leq q;\\
1, & \textrm{if} & q+1\leq i\leq p,j=j_{i-q};\\
0, & \textrm{otherwise,}\end{array}\end{cases}\end{equation}

\noindent where $1\leq j_{1}<j_{2}<\cdots<j_{p-q}\leq p$ are the
indices complementary to $1\leq i_{1}<i_{2}<\cdots<i_{q}\leq p.$
\end{singlespace}

\noindent It follows from the above definition of $\Phi_{i_{1}\cdots i_{q}}$
that

\noindent \begin{equation}
\textrm{det}\Phi_{i_{1}\cdots i_{q}}=\textrm{sgn}(\begin{array}{c}
1\cdots q\ \ q+1\cdots p\\
i_{1}\cdots i_{q}\ j_{1}\cdots j_{p-q}\end{array})\delta_{i_{1}\cdots i_{q}}.\end{equation}

\noindent We define $\tilde{u}_{i_{1}\cdots i_{q}}\in\mathcal{O}^{p}(\Omega),1\leq i_{1}<i_{2}<\cdots<i_{q}\leq p,$
by \begin{equation}
\tilde{u}_{i_{1}\cdots i_{q}}=\left(\begin{array}{c}
u_{i_{1}\cdots i_{q}}\\
v_{i_{1}\cdots i_{q}}\end{array}\right),\end{equation}

\noindent where $v_{i_{1}\cdots i_{q}}\in\mathcal{O}^{p-q}(\Omega)$
could be an arbitrary section.

\noindent By using $\tilde{u}_{i_{1}\cdots i_{q}}$ we introduce an
element of $\mathcal{O}^{p}(\Omega)$ as follows \begin{equation}
h_{i_{1}\cdots i_{q}}=\Phi_{i_{1}\cdots i_{q}}^{*}\tilde{u}_{i_{1}\cdots i_{q}}\end{equation}

\noindent for $1\leq i_{1}<i_{2}<\cdots<i_{q}\leq p.$ In the definition
(67), $\Phi_{i_{1}\cdots i_{q}}^{*}$ is the adjoint matrix of $\Phi_{i_{1}\cdots i_{q}}.$

\noindent Multiplying (67) by the matrix $\Phi_{i_{1}\cdots i_{q}},$
then the identity (65) gives

\begin{center}
$\Phi_{i_{1}\cdots i_{q}}h_{i_{1}\cdots i_{q}}=\textrm{det}\Phi_{i_{1}\cdots i_{q}}\tilde{u}_{i_{1}\cdots i_{q}}\ \ \ \ \ \ \ \ \ \ \ \ \ \ \ \ \ \ \ \ $
\par\end{center}

\begin{center}
$\ \ \ \ \ \ \ \ \ \ \ \ \ \ \ \ \ \ \ \ \ \ \ \ \ \,=\textrm{sgn}(\begin{array}{c}
1\cdots q\ \ q+1\cdots p\\
i_{1}\cdots i_{q}\ j_{1}\cdots j_{p-q}\end{array})\delta_{i_{1}\cdots i_{q}}\tilde{u}_{i_{1}\cdots i_{q}}$
\par\end{center}

\noindent By definition (64) we know

\noindent \begin{center}
$\Phi_{i_{1}\cdots i_{q}}h_{i_{1}\cdots i_{q}}=\left(\begin{array}{c}
\begin{array}{c}
\Phi h_{i_{1}\cdots i_{q}}\\*\end{array}\\
*\end{array}\right)$,
\par\end{center}

\noindent so comparing the first $q$ rows in the above equality and
(66) shows that \begin{equation}
\Phi h_{i_{1}\cdots i_{q}}=\textrm{sgn}(\begin{array}{c}
1\cdots q\ \ q+1\cdots p\\
i_{1}\cdots i_{q}\ j_{1}\cdots j_{p-q}\end{array})\delta_{i_{1}\cdots i_{q}}u_{i_{1}\cdots i_{q}}.\end{equation}

\noindent Now by using $h_{i_{1}\cdots i_{q}}\in\mathcal{O}^{p}(\Omega)$
given in (67), we set

\noindent \begin{center}
$h=\underset{i_{1}<\cdots<i_{q}}{\sum}\textrm{sgn}(\begin{array}{c}
1\cdots q\ \ q+1\cdots p\\
i_{1}\cdots i_{q}\ j_{1}\cdots j_{p-q}\end{array})h_{i_{1}\cdots i_{q}},$
\par\end{center}

\noindent then from (63) and (68) we obtain

\begin{center}
$\Phi h=\textrm{sgn}(\begin{array}{c}
1\cdots q\ \ q+1\cdots p\\
i_{1}\cdots i_{q}\ j_{1}\cdots j_{p-q}\end{array})\Phi h_{i_{1}\cdots i_{q}}\ \ \ $
\par\end{center}

\begin{center}
$\overset{(68)}{=}\underset{i_{1}<\cdots<i_{q}}{\sum}\delta_{i_{1}\cdots i_{q}}u_{i_{1}\cdots i_{q}}\overset{(63)}{=}f.\ \ \ \ \ \ \ \ $
\par\end{center}

\begin{spacing}{0.7}
\noindent which completes the proof.
\end{spacing}

\noindent \begin{flushright}
$\ $$\square$
\par\end{flushright}

\noindent $\mathbf{Remark.}$  In {[}KT71{]}, a similar condition was
used to characterize the membership for the elements of finitely
generated submodules of $A_{p}^{\oplus m}$.

\subsection*{Acknowledgment. }

I'm grateful to professor Yum-Tong Siu for proposing the question of
studying the division problem for Koszul complex and the pullback of
an exact sequence of sheaves and for his enlightening and
stimulating discussions about this subject. I would like to thank
professors Chaohao Gu and Hesheng Hu for their constant
encouragements. I also thank the department of mathematics, Harvard
university for its hospitality where most of this work was carried
out. This work is partially supported by NSFC(11171069/A010301).\\

$\,$

\noindent E-mail address: qingchunji@fudan.edu.cn

\begin{thebibliography}{V08}
\bibitem[A04]{key-1}Andersson, M. Residue currents and ideals of
holomorphic functions. Bull. Sci. math. 128 (2004), 481\textendash{}512.

\bibitem[A06]{key-2}Andersson, M. The membership problem for polynomial
ideals in terms of residue currents, Ann. Inst. Fourier 56 (2006),
101-119.

\bibitem[AG10]{key-3}Andersson, M. and Gotmark, E. Explicit representation
of membership in polynomial ideals. Math.Ann.(2010), to appear, DOI:
10.1007/s00208-010-0524-4.

\bibitem[B87]{key-4}Brownawell, W.-D. Bounds for the degrees in the
Nullstellensatz. Ann. Math. 126 (1987), 577--591.

\bibitem[D82]{key-5}Demailly, J.-P. Estimations $L^{2}$ pour l'op\'{e}ateur
$\overline{\partial}$ d'un fibr\'{e}vectoriel holomorphe
semi-positif au-dessus d'une vari\'{e}t\'{e} k\"{a}hl\'{e}rienne
compl\`{e}te. Ann. Sci. \'{E}cole Norm. Sup. (4) 15 (1982), no. 3,
457--511.

\bibitem[EL99]{key-6}Ein, L. and Lazarsfeld, R. A geometric effective
Nullstellensatz. Invent. Math. 137 (1999), no. 2, 427--448.

\bibitem[KT71]{key-7}Kelleher,J.J. and Taylor,B.A. Finitely generated
ideals in rings of analytic functions. Math. Ann. 193(1971), 225-237.

\bibitem[OT87]{key-8}Ohsawa, T. and Takegoshi, K. On the extension
of $L^{2}$ holomorphic functions. Math. Z. 195 (1987), no. 2, 197--204.

\bibitem[Siu82]{key-9}Siu, Y.-T. Complex-analyticity of harmonic
maps, vanishing and Lefschetz theorems. J. Differential Geometry.
7(1982), 55-138.

\bibitem[Siu98]{key-10}Siu, Y.-T. Invariance of plurigenera. Invent.
Math. 134 (1998), no. 3, 661--673.

\bibitem[Siu00]{key-11}Siu, Y.-T. Extension of Twisted Pluricanonical
Sections with Plurisubharmonic Weight and Invariance of
Semipositively Twisted Plurigenera for Manifolds Not Necessarily of
General Type. Complex geometry (G\"{o}ttingen, 2000), pp. 223--277.
Springer, Berlin (2002).

\bibitem[Siu04]{key-12}Siu, Y.-T. Invariance of plurigenera and torsion-freeness
of direct image sheaves of pluricanonical bundles. Finite or infinite
dimensional complex analysis and applications, 45--83, Adv. Complex
Anal. Appl., 2, Kluwer Acad. Publ., Dordrecht, 2004.

\bibitem[Siu05]{key-13}Siu, Y.-T. Multiplier ideal sheaves in complex
and algebraic geometry. Sci. China A 48(Suppl, 2005), 1--31.

\bibitem[Siu07]{key-14}Siu, Y.-T. Techniques for the analytic proof
of the finite generation of the canonical ring. Current developments
in mathematics, 2007, 177--219, Int. Press, Somerville, MA, 2009.

\bibitem[S72]{key-15}Skoda, H. Application des techniques $L^{2}$
\'{e}a la th\'{e}orie des id\'{e}aux d'une alg\`{e}bre de fonctions
holomorphes avec poids. Ann. Sci. \'{E}cole Norm. Sup. 4(5),
545--579 (1972).

\bibitem[S78]{key-16}Skoda, H. Morphismes surjectifs de fibr\'{e}s vectoriels
semi-positifs. Ann. Sci. \'{E}cole Norm. Sup. (4) 11 (1978), no. 4,
577--611.

\bibitem[V08]{key-17}Varolin, D. Division theorems and twisted complexes.
Math. Z. 259 (2008), no. 1, 1--20.

\end{thebibliography}
\end{document}